\input amstex
\documentstyle{amsppt}
%
\catcode`@=11
\redefine\output@{%
  \def\break{\penalty-\@M}\let\par\endgraf
  \ifodd\pageno\global\hoffset=105pt\else\global\hoffset=8pt\fi  
  \shipout\vbox{%
    \ifplain@
      \let\makeheadline\relax \let\makefootline\relax
    \else
      \iffirstpage@ \global\firstpage@false
        \let\rightheadline\frheadline
        \let\leftheadline\flheadline
      \else
        \ifrunheads@ 
        \else \let\makeheadline\relax
        \fi
      \fi
    \fi
    \makeheadline \pagebody \makefootline}%
  \advancepageno \ifnum\outputpenalty>-\@MM\else\dosupereject\fi
}
\catcode`\@=\active
\nopagenumbers
\def\negskp{\hskip -2pt}
\def\Alpha{\operatorname{A}}
\def\MatGrSO{\operatorname{SO}}
\def\MatGrSL{\operatorname{SL}}
\def\vtrule{\vrule height 12pt depth 6pt}
\def\vtttrule{\vrule height 12pt depth 19pt}
\def\boxit#1#2{\vcenter{\hsize=122pt\offinterlineskip\hrule
  \line{\vtttrule\hss\vtop{\hsize=120pt\centerline{#1}\vskip 5pt
  \centerline{#2}}\hss\vtttrule}\hrule}}
\def\blue#1{#1}
\catcode`#=11\def\diez{#}\catcode`#=6
\catcode`_=11\def\podcherkivanie{_}\catcode`_=8
\def\mycite#1{\cite{\blue{#1}}\immediate\special{ps:
     ShrHPSdict begin /ShrBORDERthickness 0 def}}

\def\mytag#1{%
    \tag#1}
\def\mythetag#1{\thetag{\blue{#1}}\immediate\special{ps:
     ShrHPSdict begin /ShrBORDERthickness 0 def}}
\def\myrefno#1{\no#1}
\def\myhref#1#2{\blue{#2}\immediate\special{ps:
     ShrHPSdict begin /ShrBORDERthickness 0 def}}
\def\myEarXivlink{\myhref{http://arXiv.org}{http:/\negskp/arXiv.org}}
\def\myGeoCities{\myhref{http://www.geocities.com}{GeoCities}}
\def\mytheorem#1{\csname proclaim\endcsname{Theorem #1}}

\def\mylemma#1{\csname proclaim\endcsname{Lemma #1}}

\def\mycorollary#1{\csname proclaim\endcsname{Corollary #1}}

\def\mydefinition#1{\definition{Definition #1}}
\def\mythedefinition#1{\blue{#1}\immediate\special{ps:
     ShrHPSdict begin /ShrBORDERthickness 0 def}}

\pagewidth{360pt}
\pageheight{606pt}
\topmatter
\title
On the spinor structure of the homogeneous and isotropic 
universe in closed model.
\endtitle
\author
R.~A.~Sharipov
\endauthor
\address 5 Rabochaya street, 450003 Ufa, Russia\newline
\vphantom{a}\kern 12pt Cell Phone: +7-(917)-476-93-48
\endaddress
\email \vtop to 30pt{\hsize=280pt\noindent
\myhref{mailto:r-sharipov\@mail.ru}
{r-sharipov\@mail.ru}\newline
\myhref{mailto:R\podcherkivanie Sharipov\@ic.bashedu.ru}
{R\_\hskip 1pt Sharipov\@ic.bashedu.ru}\vss}
\endemail
\urladdr
\vtop to 20pt{\hsize=280pt\noindent
\myhref{http://www.geocities.com/r-sharipov}
{http:/\negskp/www.geocities.com/r-sharipov}\newline
\myhref{http://www.freetextbooks.boom.ru/index.html}
{http:/\negskp/www.freetextbooks.boom.ru/index.html}\vss}
\endurladdr
\abstract
    The closed homogeneous and isotropic universe is considered. 
The bundles of Weyl and Dirac spinors for this universe are 
explicitly described. Some explicit formulas for the basic fields
and for the connection components in stereographic and in spherical 
coordinates are presented.
\endabstract
\endtopmatter
\loadbold
\loadeufb
\TagsOnRight
\document
\accentedsymbol\tbvartheta{\tilde{\overline{\boldsymbol\vartheta}}
\vphantom{\boldsymbol\vartheta}}

\rightheadtext{On the spinor structure \dots}
\head
1. Stereographic projections of the sphere $S^3$. 
\endhead
\parshape 21 0pt 360pt 0pt 360pt 0pt 360pt 220pt 140pt 
220pt 140pt 220pt 140pt 220pt 140pt 220pt 140pt 220pt 140pt 
220pt 140pt 220pt 140pt 220pt 140pt 220pt 140pt 220pt 140pt 
220pt 140pt 220pt 140pt 220pt 140pt 220pt 140pt 220pt 140pt 
220pt 140pt 0pt 360pt
    It is known that the closed homogeneous and isotropic universe 
is described by a manifold diffeomorphic to the Cartesian product 
of the three-dimensional sphere $S^3$ by a straight line: $M=\Bbb R
\times S^3$ (see \S\,111 and \S\,112 in \mycite{1}). The sphere $S^3$
is a manifold that can be covered with two local charts. We choose
stereographic projections from two diametrically opposite points 
(we call them North and South poles) to their equatorial hyperplane.
\vadjust{\vskip 5pt\hbox to 0pt{\kern -10pt
\includegraphics{SpStr01.eps}\hss}\vskip -5pt}The sphere
$S^3$ is naturally presented as a three-dimensional hypersurface
in the four-dimensional space $\Bbb R^4$ with the standard Euclidean
metric. Let $z$ be a point on such a sphere with the radius $R$. Then
$z/R$ is a point on the unit sphere. Let $x$ and $y$ be the
stereographic projections of $z/R$. Then the points $x$ and $y$ 
lie on some ray in the equatorial hyperplane coming out from the center 
of the sphere. It is rather easy to derive the following formulas 
relating their coordinates:
$$
\xalignat 2
	&\hskip -2em
\Vmatrix y^1\\ y^2 \\ y^3\endVmatrix=\frac{1}{|x|^2}
\,\Vmatrix x^1\\ x^2 \\ x^3\endVmatrix,
&&\Vmatrix x^1\\ x^2 \\ x^3\endVmatrix=\frac{1}{|y|^2}
\,\Vmatrix y^1\\ y^2 \\ y^3\endVmatrix.
\mytag{1.1}
\endxalignat
$$
Here $|x|^2=(x^1)^2+(x^2)^2+(x^3)^2$ and $|y|^2=(y^1)^2+(y^2)^2+(y^3)^2$.
As for the point $z$ on the sphere, its coordinates can also be expressed
through the coordinates of $x$ and $y$. Here are the formulas for the
coordinates of the point $z$:
$$
\align
&\hskip -2em
\Vmatrix z^1\\ z^2 \\ z^3\\ z^4\endVmatrix=
\frac{2\,R}{|x|^2+1}\,\Vmatrix x^1\\ x^2 \\ x^3\\ 0\endVmatrix
+\frac{|x|^2-1}{|x|^2+1}\,\Vmatrix 0\\ 0 \\ 0\\ R\endVmatrix,
\mytag{1.2}\\
&\hskip -2em
\Vmatrix z^1\\ z^2 \\ z^3\\ z^4\endVmatrix=
\frac{2\,R}{|y|^2+1}\,\Vmatrix y^1\\ y^2 \\ y^3\\ 0\endVmatrix
-\frac{|y|^2-1}{|y|^2+1}\,\Vmatrix 0\\ 0 \\ 0\\ R\endVmatrix.
\mytag{1.3}
\endalign
$$
The standard Euclidean metric in $\Bbb R^4$ is given by the formula
$$
\hskip -2em
ds^{\kern 0.5pt 2}=(dz^1)^2+(dz^2)^2+(dz^3)^2+(dz^4)^2.
\mytag{1.4}
$$
Differentiating \mythetag{1.2} and  \mythetag{1.3} and substituting
them into \mythetag{1.4}, we derive the formulas for the induced 
metric in the local coordinates $x^1,\,x^2,\,x^3$ and $y^1,\,y^2,
\,y^3$:
$$
\align
&\hskip -2em
ds^{\kern 0.5pt 2}=\frac{4\,R^2\,(dx^1)^2+4\,R^2\,(dx^2)^2
+4\,R^2\,(dx^3)^2}{\left(|x|^2+1\right)^2},
\mytag{1.5}\\
\vspace{2ex}
&\hskip -2em
ds^{\kern 0.5pt 2}=\frac{4\,R^2\,(dy^1)^2+4\,R^2\,(dy^2)^2
+4\,R^2\,(dy^3)^2}{\left(|y|^2+1\right)^2}.
\mytag{1.6}
\endalign
$$\par
     Passing from $\Bbb R^4$ to the space-time manifold $M=\Bbb R
\times S^3$, we add new coordinates $x^0$ and $y^0$ to the initial
coordinates $x^1,\,x^2,\,x^3$ and $y^1,\,y^2,\,y^3$. The transition
functions relating $x^1,\,x^2,\,x^3$ and $y^1,\,y^2,\,y^3$ are 
determined by the formulas \mythetag{1.1}. For the newly introduced 
coordinates we set
$$
\hskip -2em
x^0=y^0.
\mytag{1.7}
$$
Using the local coordinates $x^0,\,x^1,\,x^2,\,x^3$ and $y^0,\,y^1,\,y^2,
\,y^3$ and relying upon \mythetag{1.5} and \mythetag{1.6}, we introduce 
the Minkowski type metric $\bold g$ to our model $M$ of the universe:
$$
\align
&\hskip -2em
ds^{\kern 0.5pt 2}=R^2\,(dx^0)^2
-\frac{4\,R^2\,(dx^1)^2+4\,R^2\,(dx^2)^2
+4\,R^2\,(dx^3)^2}{\left(|x|^2+1\right)^2},
\mytag{1.8}\\
\vspace{2ex}
&\hskip -2em
ds^{\kern 0.5pt 2}=R^2\,(dy^0)^2
-\frac{4\,R^2\,(dy^1)^2+4\,R^2\,(dy^2)^2
+4\,R^2\,(dy^3)^2}{\left(|y|^2+1\right)^2}.
\mytag{1.9}
\endalign
$$
From now on the parameter $R$ in formulas is not a constant. 
We shall assume it to be a function of the newly introduced 
coordinates \mythetag{1.7}:
$$
\hskip -2em
R=R(x^0)=R(y^0).
\mytag{1.10}
$$
The parameter $R$ in \mythetag{1.10} is interpreted as the radius
of the sphere $S^3$ in its Euclidean realization as a hypersurface
in $\Bbb R^4$. This parameter is the only parameter describing the
evolution of the homogeneous and isotropic universe in closed model.
One can introduce the time variable $t$ through the following
formula:
$$
R\,dx^0=R\,dy^0=c\,dt\qquad\text{($c$ is the light velocity)}.
$$
Then we can write \mythetag{1.10} as $R=R(t)$. Depending on the 
function $R(t)$ we say: the universe is stable, the universe is 
expanding, or the universe is contracting. Oscillatory regimes 
are also admissible. Unlike the Newtonian mechanics, $t$ is not 
an absolute time in the universe, but the most preferable time 
variable due to the symmetry of the model.\par
      The signature of the metric $\bold g$ is $(+,-,-,-)$. Looking 
at \mythetag{1.8} and  \mythetag{1.9}, we see that the metric tensor
$\bold g$ is diagonal in the stereographic projection charts. This 
means that the coordinate frames of these two charts
$$
\xalignat 4
&\hskip -2em
\frac{\partial}{\partial x^0},
&&\frac{\partial}{\partial x^1},
&&\frac{\partial}{\partial x^2},
&&\frac{\partial}{\partial x^3},
\quad
\mytag{1.11}\\
\vspace{2ex}
&\hskip -2em
\frac{\partial}{\partial y^0},
&&\frac{\partial}{\partial y^1},
&&\frac{\partial}{\partial y^2},
&&\frac{\partial}{\partial y^3}
\quad
\mytag{1.12}
\endxalignat
$$
are orthogonal frames. However, they are not orthonormal frames.
We normalize them introducing the following two orthonormal frames:
$$
\xalignat 2
&\hskip -2em
\bold X_0=\frac{1}{R}\,\frac{\partial}{\partial x^0},
&&\bold X_1=\frac{1+|x|^2}{2\,R}\,\frac{\partial}{\partial x^1},\\
\vspace{-1.2ex}
&&&\mytag{1.13}\\
\vspace{-1.2ex}
&\hskip -2em
\bold X_2=\frac{1+|x|^2}{2\,R}\,\frac{\partial}{\partial x^2},
&&\bold X_3=\frac{1+|x|^2}{2\,R}\,\frac{\partial}{\partial x^3},\\
\vspace{3ex}
&\hskip -2em
\bold Y_0=\frac{1}{R}\,\frac{\partial}{\partial y^0},
&&\bold Y_1=\frac{1+|y|^2}{2\,R}\,\frac{\partial}{\partial y^1},\\
\vspace{-1.2ex}
&&&\mytag{1.14}\\
\vspace{-1.2ex}
&\hskip -2em
\bold Y_2=\frac{1+|y|^2}{2\,R}\,\frac{\partial}{\partial y^2},
&&\bold Y_3=\frac{1+|y|^2}{2\,R}\,\frac{\partial}{\partial y^3}.\\
\endxalignat
$$
The frames \mythetag{1.13} and  \mythetag{1.14} are orthonormal,
i\.\,e\. the metric tensor $\bold g$ and its dual metric tensor
are given by the standard Minkowski matrix
$$
\hskip -2em
g_{ij}=g^{ij}=\Vmatrix\format \l&\quad\r&\quad\r&\quad\r\\
1 &0 &0 &0\\ 0 &-1 &0 &0\\ 0 &0 &-1 &0\\ 0 &0 &0 &-1
\endVmatrix
\mytag{1.15}
$$
in both of these two frames. But unlike the frames \mythetag{1.11} 
and \mythetag{1.12}, these two frames are non-holonomic. The
vector fields $\bold X_0,\,\bold X_1,\,\bold X_2,\,\bold X_3$
and $\bold Y_0,\,\bold Y_1,\,\bold Y_2,\,\bold Y_3$ composing 
these frames do not commute with each other:
$$
\hskip -2em
[\bold X_i,\bold X_j]=\sum^3_{k=0}c^{\,k}_{ij}\,\bold X_k.
\mytag{1.16}
$$
Using \mythetag{1.13}, one can easily find the explicit formulas 
for the commutation coefficients $c^{\,k}_{ij}$ in \mythetag{1.16}. 
Most of these coefficients are zero. Below is the list of those 
coefficients $c^{\,k}_{ij}$ which are nonzero:
$$
\xalignat 2
&\hskip -2em
\kern 1em c^{\,1}_{01}=-c^{\,1}_{10}=c^{\,2}_{02}=-c^{\,2}_{20}=
c^{\,3}_{03}=-c^{\,3}_{30}=-\frac{R'}{R^2},\kern -7cm
\quad\\
\vspace{1ex}
&\hskip -2em
c^{\,1}_{12}=-c^{\,1}_{21}=-\frac{(x^2)}{R^2},
&&c^{\,2}_{12}=-c^{\,2}_{21}=\frac{(x^1)}{R^2},
\quad\\
\vspace{-1ex}
&&&\mytag{1.17}\\
\vspace{-1ex}
&\hskip -2em
c^{\,1}_{13}=-c^{\,1}_{31}=-\frac{(x^3)}{R^2},
&&c^{\,3}_{13}=-c^{\,3}_{31}=\frac{(x^1)}{R^2},
\quad\\
\vspace{1ex}
&\hskip -2em
c^{\,2}_{23}=-c^{\,2}_{32}=-\frac{(x^3)}{R^2},
&&c^{\,3}_{23}=-c^{\,3}_{32}=\frac{(x^2)}{R^2}.
\quad
\endxalignat
$$
Here $R'$ is the derivative of the function \mythetag{1.10}.
The vector fields of the second frame $\bold Y_0,\,\bold Y_1,
\,\bold Y_2,\,\bold Y_3$ obey the commutation relationships
similar to \mythetag{1.16}:
$$
\hskip -2em
[\bold Y_i,\bold Y_j]=\sum^3_{k=0}c^{\,k}_{ij}\,\bold Y_k.
\mytag{1.18}
$$
Below is the list of all nonzero commutation coefficients 
$c^{\,k}_{ij}$ for \mythetag{1.18}:
$$
\xalignat 2
&\hskip -2em
\kern 1em c^{\,1}_{01}=-c^{\,1}_{10}=c^{\,2}_{02}=-c^{\,2}_{20}=
c^{\,3}_{03}=-c^{\,3}_{30}=-\frac{R'}{R^2},\kern -7cm
\quad\\
\vspace{1ex}
&\hskip -2em
c^{\,1}_{12}=-c^{\,1}_{21}=-\frac{(y^2)}{R^2},
&&c^{\,2}_{12}=-c^{\,2}_{21}=\frac{(y^1)}{R^2},
\quad\\
\vspace{-1ex}
&&&\mytag{1.19}\\
\vspace{-1ex}
&\hskip -2em
c^{\,1}_{13}=-c^{\,1}_{31}=-\frac{(y^3)}{R^2},
&&c^{\,3}_{13}=-c^{\,3}_{31}=\frac{(y^1)}{R^2},
\quad\\
\vspace{1ex}
&\hskip -2em
c^{\,2}_{23}=-c^{\,2}_{32}=-\frac{(y^3)}{R^2},
&&c^{\,3}_{23}=-c^{\,3}_{32}=\frac{(y^2)}{R^2}.
\quad
\endxalignat
$$\par
     The formula \mythetag{1.1} complemented with the formula 
\mythetag{1.7} determines the transition functions for two 
overlapping local charts $x^0,\,x^1,\,x^2,\,x^3$ and $y^0,\,y^1,
\,y^2,\,y^3$. Differentiating these transition functions, we get 
the transition matrices relating the holonomic frames \mythetag{1.11} 
and \mythetag{1.12}:
$$
\xalignat 2
&\hskip -2em
\frac{\partial}{\partial x^i}=
\sum^3_{j=0}\frac{\partial y^j}{\partial x^i}\,\frac{\partial}
{\partial y^j},
&&\frac{\partial}{\partial y^i}=
\sum^3_{j=0}\frac{\partial x^j}{\partial y^i}\,\frac{\partial}
{\partial x^j}.\quad
\mytag{1.20}
\endxalignat
$$
Then, applying \mythetag{1.13} and \mythetag{1.14} to \mythetag{1.20}, 
we derive the formulas relating the non-holonomic frames $\bold X_0,
\,\bold X_1,\,\bold X_2,\,\bold X_3$ and $\bold Y_0,\,\bold Y_1,
\,\bold Y_2,\,\bold Y_3$:
$$
\xalignat 2
&\hskip -2em
\bold X_i=\sum^3_{j=0}T^j_i\,\bold Y_j,
&&\bold Y_i=
\sum^3_{j=0}S^j_i\,\bold X_j.\quad
\mytag{1.21}
\endxalignat
$$
Here is the explicit formula for the matrix $S$ in \mythetag{1.21}:
$$
\hskip -2em
S=\Vmatrix
1 & 0 & 0 & 0\\
\vspace{3ex}
0 & \dfrac{|y|^2-2\,(y^1)^2}{|y|^2}
& \dfrac{-2\,(y^1)\,(y^2)}{|y|^2}
& \dfrac{-2\,(y^1)\,(y^3)}{|y|^2}\\
\vspace{3ex}
0 & \dfrac{-2\,(y^1)\,(y^2)}{|y|^2}
& \dfrac{|y|^2-2\,(y^2)^2}{|y|^2} 
& \dfrac{-2\,(y^2)\,(y^3)}{|y|^2}\\
\vspace{3ex}
0 & \dfrac{-2\,(y^1)\,(y^3)}{|y|^2}
& \dfrac{-2\,(y^2)\,(y^3)}{|y|^2} 
& \dfrac{|y|^2-2\,(y^3)^2}{|y|^2}\\
\endVmatrix.
\mytag{1.22}
$$
The matrix $T$ is the inverse matrix for $S$, i\.\,e\. $T=S^{-1}$.
Its components can be expressed in terms of the coordinates 
$x^0,\,x^1,\,x^2,\,x^3$:
$$
\hskip -2em
T=\Vmatrix
1 & 0 & 0 & 0\\
\vspace{3ex}
0 & \dfrac{|x|^2-2\,(x^1)^2}{|x|^2}
& \dfrac{-2\,(x^1)\,(x^2)}{|x|^2}
& \dfrac{-2\,(x^1)\,(x^3)}{|x|^2}\\
\vspace{3ex}
0 & \dfrac{-2\,(x^1)\,(x^2)}{|x|^2}
& \dfrac{|x|^2-2\,(x^2)^2}{|x|^2} 
& \dfrac{-2\,(x^2)\,(x^3)}{|x|^2}\\
\vspace{3ex}
0 & \dfrac{-2\,(x^1)\,(x^3)}{|x|^2}
& \dfrac{-2\,(x^2)\,(x^3)}{|x|^2} 
& \dfrac{|x|^2-2\,(x^3)^2}{|x|^2}\\
\endVmatrix.
\mytag{1.23}
$$
In general case the matrices $S$ and $T$ are different. However, 
it is the feature of our particular charts that the matrices 
\mythetag{1.22} and \mythetag{1.23} do coincide:
$$
\xalignat 2
&\hskip -2em
S=T, &&S^2=T^2=1.
\mytag{1.24}
\endxalignat
$$\par
     The frames \mythetag{1.13} and \mythetag{1.14} are two orthonormal
frames of the spherical universe $M=\Bbb R\times S^3$. The vector $\bold
X_0=\bold Y_0$ is a time-like vector directed to the future. Therefore 
the frames \mythetag{1.13} and \mythetag{1.14} are called {\it positively 
polarized} (see definition in \S\,5 of \mycite{2}). However, their
orientations are different. Indeed, by means of direct calculations for
$S$ and $T$ we find that 
$$
\hskip -2em
\det S=\det T=-1.
\mytag{1.25}
$$
The formula \mythetag{1.25} is concordant with \mythetag{1.24}. It means
that if we take \mythetag{1.13} for a right oriented frame in $M$, then
\mythetag{1.14} is a left oriented frame. In the theory of Weyl spinors
only positively polarized right orthonormal frames are admissible (see
definition~5.2 in \mycite{2}). For this reason we introduce the following
auxiliary frame:
$$
\pagebreak
\xalignat 4
&\hskip -2em
\tilde\bold Y_0=\bold Y_0,
&&\tilde\bold Y_1=-\bold Y_1,
&&\tilde\bold Y_2=-\bold Y_2,
&&\tilde\bold Y_3=-\bold Y_3.
\qquad
\mytag{1.26}
\endxalignat
$$
Like $\bold X_0,\,\bold X_1,\,\bold X_2,\,\bold X_3$, the frame
\mythetag{1.26} is a positively polarized right orthonormal frame 
in $M$. It is related to $\bold X_0,\,\bold X_1,\,\bold X_2,
\,\bold X_3$ as follows:
$$
\xalignat 2
&\hskip -2em
\bold X_i=\sum^3_{j=0}\tilde T^j_i\,\tilde\bold Y_j,
&&\tilde\bold Y_i=
\sum^3_{j=0}\tilde S^j_i\,\bold X_j.\quad
\mytag{1.27}
\endxalignat
$$
Here are the explicit formulas for the matrices $\tilde S$ and
$\tilde T$ in \mythetag{1.27}:
$$
\align
&\hskip -2em
\tilde S=\Vmatrix
1 & 0 & 0 & 0\\
\vspace{3ex}
0 & \dfrac{2\,(y^1)^2-|y|^2}{|y|^2}
& \dfrac{2\,(y^1)\,(y^2)}{|y|^2}
& \dfrac{2\,(y^1)\,(y^3)}{|y|^2}\\
\vspace{3ex}
0 & \dfrac{2\,(y^1)\,(y^2)}{|y|^2}
& \dfrac{2\,(y^2)^2-|y|^2}{|y|^2} 
& \dfrac{2\,(y^2)\,(y^3)}{|y|^2}\\
\vspace{3ex}
0 & \dfrac{2\,(y^1)\,(y^3)}{|y|^2}
& \dfrac{2\,(y^2)\,(y^3)}{|y|^2} 
& \dfrac{2\,(y^3)^2-|y|^2}{|y|^2}
\endVmatrix,
\mytag{1.28}\\
\vspace{2ex}
&\hskip -2em
\tilde T=\Vmatrix
1 & 0 & 0 & 0\\
\vspace{3ex}
0 & \dfrac{|x|^2-2\,(x^1)^2}{|x|^2}
& \dfrac{-2\,(x^1)\,(x^2)}{|x|^2}
& \dfrac{-2\,(x^1)\,(x^3)}{|x|^2}\\
\vspace{3ex}
0 & \dfrac{-2\,(x^1)\,(x^2)}{|x|^2}
& \dfrac{|x|^2-2\,(x^2)^2}{|x|^2} 
& \dfrac{-2\,(x^2)\,(x^3)}{|x|^2}\\
\vspace{3ex}
0 & \dfrac{-2\,(x^1)\,(x^3)}{|x|^2}
& \dfrac{-2\,(x^2)\,(x^3)}{|x|^2} 
& \dfrac{|x|^2-2\,(x^3)^2}{|x|^2}
\endVmatrix.
\mytag{1.29}
\endalign
$$
Like in \mythetag{1.24}, for the matrices \mythetag{1.28} and 
\mythetag{1.29} we have the relationships
$$
\xalignat 2
&\hskip -2em
\tilde S=\tilde T, &&\tilde S^2=\tilde T^2=1.
\mytag{1.30}
\endxalignat
$$
However, instead of \mythetag{1.25}, their determinants now are equal 
to the unity:
$$
\hskip -2em
\det\tilde S=\det\tilde T=1.
\mytag{1.31}
$$
Being transition matrices that relate two positively polarized 
right orthonormal frames, the matrices $\tilde S$ and $\det\tilde T$ 
belong to the special orthochronous Lorentz group 
$\MatGrSO^+(1,3,\Bbb R)$. Our next step is to construct the bundle 
$SM$ of Weyl spinors for the universe $M=\Bbb R\times S^3$. We use 
the frames $\bold X_0,\,\bold X_1,\,\bold X_2,\,\bold X_3$ and 
$\tilde\bold Y_0,\,\tilde\bold Y_1,\,\tilde\bold Y_2,\,\tilde
\bold Y_3$ and their transition matrices $\tilde S$ and $\tilde T$
for this purpose.\par
\head
2. Constructing the bundle of Weyl spinors.
\endhead
    The bundle of Weyl spinors $SM$ is a two-dimensional complex
vector-bundle related in some special way to the tangent bundle
$TM$ (see definition~5.2 in \mycite{2}). According to this
definition, each positively polarized right orthonormal frame 
of the tangent bundle $TM$ should be associated with some frame 
of the spinor bundle in such a way that the transition matrices 
of the associated spinor frames would belong to the group 
$\MatGrSL(2,\Bbb C)$ and would be linked to the transition 
matrices of the tangent frames by means of the group homomorphism
$$
\hskip -2em
\phi\!:\,\MatGrSL(2,\Bbb C)\to\MatGrSO^+(1,3,\Bbb R).
\mytag{2.1}
$$
In our particular case we have two positively polarized right
orthonormal frames \mythetag{1.13} and \mythetag{1.26} with the
transition matrix \mythetag{1.28} relating them. In order to prove 
the existence of the spinor bundle $SM$ in the case of the spherical
universe $M=\Bbb R\times S^3$ we need to find 
a matrix $\tilde{\goth S}\in\MatGrSL(2,\Bbb C)$ such that
$\phi(\tilde{\goth S})=\tilde S$. The components of the matrix 
$\tilde{\goth S}$ should be smooth functions in the intersection
of the domains of two frames \mythetag{1.13} and \mythetag{1.26},
i\.\,e\. they should be smooth functions on the whole sphere 
$S^3$ except for the poles. As for the homomorphism \mythetag{2.1}, 
this homomorphism is given by the explicit formulas \thetag{1.2}, 
\thetag{1.3}, \thetag{1.4}, and \thetag{1.5} in paper \mycite{3}.
Here are these explicit formulas expressing $\tilde S
=\phi(\tilde{\goth S})$ through $\tilde{\goth S}$:
$$
\allowdisplaybreaks
\gather
\hskip -6em
\gathered
\tilde S^0_0=\frac{\overline{\tilde\goth S^1_1}\,\tilde\goth S^1_1
+\overline{\tilde\goth S^1_2}\,\tilde\goth S^1_2+\overline{\tilde
\goth S^2_1}\,\tilde\goth S^2_1+\overline{\tilde\goth S^2_2}\,
\tilde\goth S^2_2}{2},\\
\tilde S^0_1=\frac{\overline{\tilde\goth S^1_1}\,\tilde\goth S^1_2
+\overline{\tilde\goth S^1_2}\,\tilde\goth S^1_1+\overline{\tilde
\goth S^2_1}\,\tilde\goth S^2_2+\overline{\tilde\goth S^2_2}\,
\tilde\goth S^2_1}{2},\\
\tilde S^0_2=\frac{\overline{\tilde\goth S^1_2}\,\tilde\goth S^1_1
-\overline{\tilde\goth S^1_1}\,\tilde\goth S^1_2+\overline{\tilde
\goth S^2_2}\,\tilde\goth S^2_1-\overline{\tilde\goth S^2_1}\,
\tilde\goth S^2_2}{2\,i},\\
\tilde S^0_3=\frac{\overline{\tilde\goth S^1_1}\,\tilde\goth S^1_1
-\overline{\tilde\goth S^1_2}\,\tilde\goth S^1_2+\overline{\tilde
\goth S^2_1}\,\tilde\goth S^2_1-\overline{\tilde\goth S^2_2}\,
\tilde\goth S^2_2}{2},
\endgathered
\mytag{2.2}\\
\vspace{2ex}
\hskip 2em
\gathered
\tilde S^1_0=\frac{\overline{\tilde\goth S^2_1}\,\tilde\goth S^1_1
+\overline{\tilde\goth S^1_1}\,\tilde\goth S^2_1+\overline{\tilde
\goth S^2_2}\,\tilde\goth S^1_2+\overline{\tilde\goth S^1_2}\,
\tilde\goth S^2_2}{2},\\
\tilde S^1_1=\frac{\overline{\tilde\goth S^2_1}\,\tilde\goth S^1_2
+\overline{\tilde\goth S^1_2}\,\tilde\goth S^2_1+\overline{\tilde
\goth S^2_2}\,\tilde\goth S^1_1+\overline{\tilde\goth S^1_1}\,
\tilde\goth S^2_2}{2},\\
\tilde S^1_2=\frac{\overline{\tilde\goth S^1_2}\,\tilde\goth S^2_1
-\overline{\tilde\goth S^2_1}\,\tilde\goth S^1_2+\overline{\tilde
\goth S^2_2}\,\tilde\goth S^1_1-\overline{\tilde\goth S^1_1}\,
\tilde\goth S^2_2}{2\,i},\\
\tilde S^1_3=\frac{\overline{\tilde\goth S^2_1}\,\tilde\goth S^1_1
+\overline{\tilde\goth S^1_1}\,\tilde\goth S^2_1-\overline{\tilde
\goth S^2_2}\,\tilde\goth S^1_2-\overline{\tilde\goth S^1_2}\,
\tilde\goth S^2_2}{2},
\endgathered
\mytag{2.3}\\
\vspace{2ex}
\hskip -6em
\gathered
\tilde S^2_0=\frac{\overline{\tilde\goth S^1_1}\,\tilde\goth S^2_1
-\overline{\tilde\goth S^2_1}\,\tilde\goth S^1_1+\overline{\tilde
\goth S^1_2}\,\tilde\goth S^2_2-\overline{\tilde\goth S^2_2}\,
\tilde\goth S^1_2}{2\,i},\\
\tilde S^2_1=\frac{\overline{\tilde\goth S^1_2}\,\tilde\goth S^2_1
-\overline{\tilde\goth S^2_1}\,\tilde\goth S^1_2+\overline{\tilde
\goth S^1_1}\,\tilde\goth S^2_2-\overline{\tilde\goth S^2_2}\,
\tilde\goth S^1_1}{2\,i},\\
\tilde S^2_2=\frac{\overline{\tilde\goth S^2_2}\,\tilde\goth S^1_1
+\overline{\tilde\goth S^1_1}\,\tilde\goth S^2_2-\overline{\tilde
\goth S^2_1}\,\tilde\goth S^1_2-\overline{\tilde\goth S^1_2}\,
\tilde\goth S^2_1}{2},\\
\tilde S^2_3=\frac{\overline{\tilde\goth S^1_1}\,\tilde\goth S^2_1
-\overline{\tilde\goth S^2_1}\,\tilde\goth S^1_1+\overline{\tilde
\goth S^2_2}\,\tilde\goth S^1_2-\overline{\tilde\goth S^1_2}\,
\tilde\goth S^2_2}{2\,i},
\endgathered
\mytag{2.4}\\
\vspace{2ex}
\hskip 2em
\gathered
\tilde S^3_0=\frac{\overline{\tilde\goth S^1_1}\,\tilde\goth S^1_1
+\overline{\tilde\goth S^1_2}\,\tilde\goth S^1_2-\overline{\tilde
\goth S^2_1}\,\tilde\goth S^2_1-\overline{\tilde\goth S^2_2}\,
\tilde\goth S^2_2}{2},\\
\tilde S^3_1=\frac{\overline{\tilde\goth S^1_1}\,\tilde\goth S^1_2
+\overline{\tilde\goth S^1_2}\,\tilde\goth S^1_1-\overline{\tilde
\goth S^2_1}\,\tilde\goth S^2_2-\overline{\tilde\goth S^2_2}\,
\tilde\goth S^2_1}{2},\\
\tilde S^3_2=\frac{\overline{\tilde\goth S^1_2}\,\tilde\goth S^1_1
-\overline{\tilde\goth S^1_1}\,\tilde\goth S^1_2+\overline{\tilde
\goth S^2_1}\,\tilde\goth S^2_2-\overline{\tilde\goth S^2_2}\,
\tilde\goth S^2_1}{2\,i},\\
\tilde S^3_3=\frac{\overline{\tilde\goth S^1_1}\,\tilde\goth S^1_1
+\overline{\tilde\goth S^2_2}\,\tilde\goth S^2_2-\overline{\tilde
\goth S^2_1}\,\tilde\goth S^2_1-\overline{\tilde\goth S^1_2}\,
\tilde\goth S^1_2}{2}.
\endgathered
\mytag{2.5}
\endgather
$$
The components of the matrix $\tilde S$ are known. Therefore,
the formulas \mythetag{2.2}, \mythetag{2.3}, \mythetag{2.4}, and
\mythetag{2.5}, are understood as the equations for the components 
of a complex $2\times 2$ matrix $\tilde\goth S$. As appears, these
equations can be solved explicitly:
$$
\hskip -2em
\tilde\goth S=\frac{1}{|y|}
\Vmatrix
i\,y^3 & i\,y^1+y^2\\
\vspace{3ex}
i\,y^1-y^2 & -i\,y^3
\endVmatrix.
\mytag{2.6}
$$
It is easy to see that $\det\tilde\goth S=1$, which means that 
$\tilde\goth S\in\MatGrSL(2,\Bbb C)$. The matrix \mythetag{2.6}
satisfying the equations \mythetag{2.2}, \mythetag{2.3}, 
\mythetag{2.4}, \mythetag{2.5} and belonging to $\MatGrSL(2,\Bbb C)$
is unique up to the change of sign: $\tilde\goth S\to -\tilde\goth S$.
\par
     Let's denote by $\tilde\goth T$ the inverse matrix for the 
matrix \mythetag{2.6}, i\.\,e\. let $\tilde\goth T=\tilde\goth S^{-1}$.
By means of the direct calculations, applying \mythetag{1.1},
we find
$$
\hskip -2em
\tilde\goth T=\frac{-1}{|x|}
\Vmatrix
i\,x^3 & i\,x^1+x^2\\
\vspace{3ex}
i\,x^1-x^2 & -i\,x^3
\endVmatrix.
\mytag{2.7}
$$
When expressed back through $y^0,\,y^1,\,y^2,\,y^3$, the matrix 
\mythetag{2.7} coincides with $-\tilde\goth S$. This means that
we have the relationships
$$
\xalignat 2
&\hskip -2em
\tilde\goth S=-\tilde\goth T, &&\tilde\goth S^2=\tilde\goth T^2=-1.
\mytag{2.8}
\endxalignat
$$
The determinants of the matrices \mythetag{2.6} and 
\mythetag{2.7} are equal to the unity:
$$
\hskip -2em
\det\tilde S=\det\tilde T=1.
\mytag{2.9}
$$
The relationships \mythetag{2.8} and \mythetag{2.9} are concordant
with the relationships \mythetag{1.30} and \mythetag{1.31} since 
$\phi(-\goth S)=\phi(\goth S)$ for the homomorphism 
\mythetag{2.1}.\par
     The mutually inverse matrices \mythetag{2.6} and \mythetag{2.7} 
are postulated to be the transition matrices for the spinor frames 
$\boldsymbol\Psi_1,\,\boldsymbol\Psi_2$ and $\tilde{\boldsymbol\Phi}_1,
\,\tilde{\boldsymbol\Phi}_2$:
$$
\xalignat 2
&\hskip -2em
\boldsymbol\Psi_i=\sum^2_{j=1}\tilde\goth T^j_i\,
\tilde{\boldsymbol\Phi}_j,
&&\tilde{\boldsymbol\Phi}_i=
\sum^2_{j=1}\tilde\goth S^j_i\,\boldsymbol\Psi_j.\quad
\mytag{2.10}
\endxalignat
$$
Thus, having found the matrices \mythetag{2.6} and \mythetag{2.7}
and having equipped the local charts $x^0,\,x^1,\,x^2,\,x^3$ and
$y^0,\,y^1,\,y^2,\,y^3$ \pagebreak 
with the spinor frames $\boldsymbol\Psi_1,\,\boldsymbol\Psi_2$ and 
$\tilde{\boldsymbol\Phi}_1,
\,\tilde{\boldsymbol\Phi}_2$ related to each other by means of the 
formulas \mythetag{2.10}, we have constructed the spinor bundle 
$SM$ over the space-time manifold $M=\Bbb R\times S^3$.
\head
3. Basic fields of the bundle of Weyl spinors.
\endhead
     The spinor bundle of Weyl spinors $SM$ over any four-dimensional
space-time manifold $M$ is equipped with two special spin-tensorial
fields. These basic spin-tensorial fields are presented in the following 
table:
$$
\vcenter{\hsize 10cm
\offinterlineskip\settabs\+\indent
\vtrule
\hskip 1.2cm &\vtrule 
\hskip 5.2cm &\vtrule 
\hskip 2.8cm &\vtrule 
\cr\hrule 
\+\vtrule
\hfill\,Symbol\hfill&\vtrule
\hfill Name\hfill &\vtrule
\hfill Spin-tensorial\hfill &\vtrule\cr
\vskip -0.2cm
\+\vtrule
\hfill &\vtrule
\hfill \hfill&\vtrule
\hfill type\hfill&\vtrule\cr\hrule
\+\vtrule
\hfill $\bold d$\hfill&\vtrule
\hfill Skew-symmetric metric tensor\hfill&\vtrule
\hfill $(0,2|0,0|0,0)$\hfill&\vtrule\cr\hrule
\+\vtrule
\hfill$\bold G$\hfill&\vtrule
\hfill Infeld-van der Waerden field\hfill&\vtrule
\hfill $(1,0|1,0|0,1)$\hfill&\vtrule\cr\hrule
}\quad
\mytag{3.1}
$$
The spin-tensorial type in the table \mythetag{3.1} specifies
the number of indices in coordinate representation of fields. 
The first two numbers are the numbers of upper and lower spinor 
indices, the second two numbers are the numbers of upper and 
lower conjugate spinor indices, and the last two numbers are 
the numbers of upper and lower tensorial indices (they are also 
called spacial indices).\par
     Now let's return to our special case $M=\Bbb R\times S^3$.
The spinor frames $\boldsymbol\Psi_1,\,\boldsymbol\Psi_2$ and 
$\tilde{\boldsymbol\Phi}_1,\,\tilde{\boldsymbol\Phi}_2$ 
considered in section 2 are canonically associated with positively
polarized right orthonormal frames in $TM$. For this reason
they are orthonormal frames by definition, i\.\,e\. the skew
symmetric metric tensor $\bold d$ is given by the matrix
$$
\hskip -2em
d_{ij}=\Vmatrix 0 & 1\\ 
\vspace{1ex} -1 & 0\endVmatrix
\mytag{3.2}
$$
in both of these frames. The indices $i$ and $j$ in \mythetag{3.2} 
are spinor indices. Therefore, the canonical presentation 
\mythetag{3.2} of the metric tensor $\bold d$ depends on the choice 
of a spinor frame, but it is not sensitive to the choice of a
tangent frame.\par
     Unlike $\bold d$, the Infeld-van der Waerden field $\bold G$
has one lower spacial index in its coordinate presentation. Therefore,
its coordinate presentation depends on the choice of two frames in
$TM$ and in $SM$. According to the definition of the spinor bundle
$SM$ (see definition~5.2 in \mycite{2}), each positively polarized 
right orthonormal frame of $TM$ has its associated orthonormal frame
in $SM$. We visualize this frame association through the following
diagram:
$$
\hskip -2em
\aligned
&\boxit{\lower 5pt\hbox{Orthonormal frames}}{}\to
\boxit{Positively polarized}{right orthonormal frames}
\endaligned\quad
\mytag{3.3}
$$
In each canonically associated pair of frames the components of
the Infeld-van der Waerden field $G^{i\kern 0.5pt\bar i}_q$ are
presented by the Pauli matrices:
$$
\xalignat 2
&\hskip -2em
G^{i\kern 0.5pt\bar i}_0=\Vmatrix 1 & 0\\ 0 & 1\endVmatrix,
&&G^{i\kern 0.5pt\bar i}_2=\Vmatrix 0 & -i\\ i & 0\endVmatrix,\\
\vspace{-1.4ex}
&&&\mytag{3.4}\\
\vspace{-1.4ex}
&\hskip -2em
G^{i\kern 0.5pt\bar i}_1=\Vmatrix 0 & 1\\ 1 & 0\endVmatrix,
&&G^{i\kern 0.5pt\bar i}_3=\Vmatrix 1 & 0\\ 0 & -1\endVmatrix.\\
\endxalignat
$$
The lower spacial index $q=0,\,1,\,2,\,3$ enumerates the matrices
in \mythetag{3.4}. The spinor index $i$ and the conjugate spinor
index $\bar i$ determine the position of the component $G^{i\kern 
0.5pt \bar i}_q$ within one of these matrices. In section 2 above 
we have constructed the following two pairs of associated frames 
of $SM$ and $TM$:
$$
\align
\hskip -2em
\boldsymbol\Psi_1,\,\boldsymbol\Psi_2\qquad&\to\qquad
\bold X_0,\,\bold X_1,\,\bold X_2,\,\bold X_3,
\mytag{3.5}\\
\hskip -2em
\tilde{\boldsymbol\Phi}_1,\,\tilde{\boldsymbol\Phi}_2
\qquad&\to\qquad
\tilde\bold Y_0,\,\tilde\bold Y_1,\,\tilde\bold Y_2,
\,\tilde\bold Y_3.
\mytag{3.6}
\endalign
$$
Once the spinor bundle $SM$ is constructed, the choice of associated 
frame pairs is not obligatory, e\.\,g\. we can choose the frames
$$
\hskip -2em
\tilde{\boldsymbol\Phi}_1,\,\tilde{\boldsymbol\Phi}_2
\qquad-\qquad
\bold Y_0,\,\bold Y_1,\,\bold Y_2,\,\bold Y_3.
\mytag{3.7}
$$
These frames are not canonically associated. However, using 
\mythetag{1.26} and \mythetag{3.4}, we easily calculate the 
Infeld-van der Waerden symbols $G^{i\kern 0.5pt\bar i}_q$ in 
this frame pair:
$$
\xalignat 2
&G^{i\kern 0.5pt\bar i}_0=\Vmatrix 1 & 0\\ 0 & 1\endVmatrix,
&&G^{i\kern 0.5pt\bar i}_2=-\Vmatrix 0 & -i\\ i & 0\endVmatrix,\\
\vspace{1ex}
&G^{i\kern 0.5pt\bar i}_1=-\Vmatrix 0 & 1\\ 1 & 0\endVmatrix,
&&G^{i\kern 0.5pt\bar i}_3=-\Vmatrix 1 & 0\\ 0 & -1\endVmatrix.\\
\endxalignat
$$
Similarly, one can combine any frame of $SM$ with any frame of $TM$
into a frame pair and then calculate the components of $\bold G$ 
in such a non-canonical frame pair.
\head
4. The bundle of Dirac spinors and its basic fields. 
\endhead
     The bundle of Dirac spinors $DM$ is constructed as the direct 
sum of the bundle of Weyl spinors $SM$ and its Hermitian conjugate 
bundle $S^{\sssize\dagger}\!M$:
$$
\hskip -2em
DM=SM\oplus S^{\sssize\dagger}\!M.
\mytag{4.1}
$$
It does exist provided $SM$ does. The bundles $SM$ and
$S^{\sssize\dagger}\!M$ are called chiral and antichiral components 
of the expansion \mythetag{4.1}. The Dirac bundle $DM$ has more
basic spin-tensorial fields as compared to $SM$:
$$
\vcenter{\hsize 10cm
\offinterlineskip\settabs\+\indent
\vtrule
\hskip 1.2cm &\vtrule 
\hskip 5.2cm &\vtrule 
\hskip 2.8cm &\vtrule 
\cr\hrule 
\+\vtrule
\hfill\,Symbol\hfill&\vtrule
\hfill Name\hfill &\vtrule
\hfill Spin-tensorial\hfill &\vtrule\cr
\vskip -0.2cm
\+\vtrule
\hfill &\vtrule
\hfill \hfill&\vtrule
\hfill type\hfill&\vtrule\cr\hrule
\+\vtrule
\hfill $\bold d$\hfill&\vtrule
\hfill Skew-symmetric metric tensor\hfill&\vtrule
\hfill $(0,2|0,0|0,0)$\hfill&\vtrule\cr\hrule
\+\vtrule
\hfill$\bold H$\hfill&\vtrule
\hfill Chirality operator\hfill&\vtrule
\hfill $(1,1|0,0|0,0)$\hfill&\vtrule\cr\hrule
\+\vtrule
\hfill$\bold D$\hfill&\vtrule
\hfill Dirac form\hfill&\vtrule
\hfill $(0,1|0,1|0,0)$\hfill&\vtrule\cr\hrule
\+\vtrule
\hfill$\boldsymbol\gamma$\hfill&\vtrule
\hfill Dirac $\gamma$-field\hfill&\vtrule
\hfill $(1,1|0,0|0,1)$\hfill&\vtrule\cr\hrule
}\quad
\mytag{4.2}
$$
Let $\boldsymbol\Psi_1,\,\boldsymbol\Psi_2$ be an orthonormal frame
of the bundle $SM$ and let $\overline{\boldsymbol\vartheta}
\vphantom{\boldsymbol\vartheta}^{\,1},\overline{\boldsymbol\vartheta}
\vphantom{\boldsymbol\vartheta}^{\,2}$ be its conjugate frame in
$S^{\sssize\dagger}\!M$. At each point $p$ off the space-time manifold
$\overline{\boldsymbol\vartheta}\vphantom{\boldsymbol\vartheta}^{\,1}$
and $\overline{\boldsymbol\vartheta}\vphantom{\boldsymbol
\vartheta}^{\,2}$ are semilinear functionals in the fiber $S_p(M)$ 
such that
$$
\overline{\boldsymbol\vartheta}\vphantom{\boldsymbol\vartheta}^{\,i}
(\boldsymbol\Psi_1)=\delta^i_j,
$$
where $\delta^i_j$ is the Kronecker delta-symbol. Denoting 
$\boldsymbol\Psi_3=\overline{\boldsymbol\vartheta}
\vphantom{\boldsymbol\vartheta}^{\,1}$ and $\boldsymbol\Psi_4
=\overline{\boldsymbol\vartheta}\vphantom{\boldsymbol
\vartheta}^{\,2}$, we get a frame $\boldsymbol\Psi_1,\,\boldsymbol\Psi_2,
\,\boldsymbol\Psi_3,\,\boldsymbol\Psi_4$ of the Dirac bundle 
\mythetag{4.1} (see more details in \mycite{3}). Such a frame of $DM$
is called a canonically orthonormal chiral frame. There are four
types of special frames in $DM$. Each special frame is canonically 
associated with some frame in tangent bundle $TM$ according to the
following diagram:
$$
\hskip -2em
\aligned
&\boxit{Canonically orthonormal}{chiral frames}\to
\boxit{Positively polarized}{right orthonormal frames}\\
&\boxit{$P$-reverse}{antichiral frames}\to
\boxit{Positively polarized}{left orthonormal frames}\\
&\boxit{$T$-reverse}{antichiral frames}\to
\boxit{Negatively polarized}{right orthonormal frames}\\
&\boxit{$PT$-reverse}{chiral frames}\to
\boxit{Negatively polarized}{left orthonormal frames}
\endaligned
\mytag{4.3}
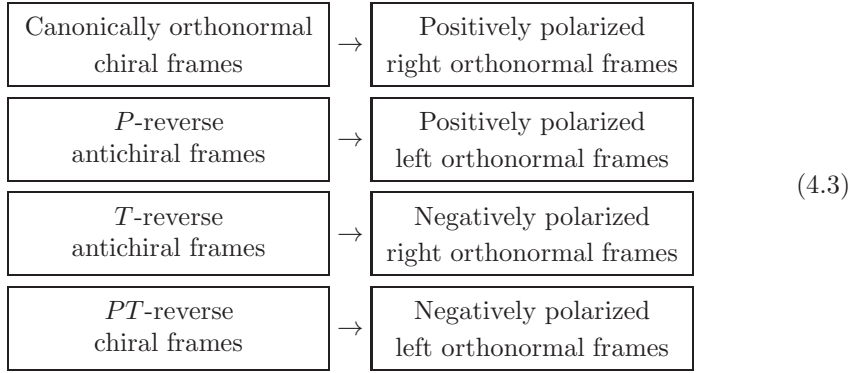
$$
\mydefinition{4.1} A frame $\boldsymbol\Psi_1,\,\boldsymbol\Psi_2,
\,\boldsymbol\Psi_3,\,\boldsymbol\Psi_4$ of the Dirac bundle is 
called an orthonormal frame if the metric tensor $\bold d$ is given
by the matrix 
$$
\hskip -2em
d_{ij}=\Vmatrix 0 & 1 & 0 & 0\\-1 & 0 & 0 & 0\\
0 & 0 & 0 & -1\\0 & 0 & 1 & 0\endVmatrix
\mytag{4.4}
$$
in this frame. This frame is called an anti-orthonormal frame if 
the metric tensor $\bold d$ is given by the opposite matrix in
this frame:
$$
\hskip -2em
d_{ij}=-\Vmatrix 0 & 1 & 0 & 0\\-1 & 0 & 0 & 0\\
0 & 0 & 0 & -1\\0 & 0 & 1 & 0\endVmatrix.
\mytag{4.5}
$$
\enddefinition
\mydefinition{4.2} A frame $\boldsymbol\Psi_1,\,\boldsymbol\Psi_2,
\,\boldsymbol\Psi_3,\,\boldsymbol\Psi_4$ of the Dirac bundle is 
called a chiral frame if the chirality operator $\bold H$ is given
by the matrix 
$$
H^i_j=\Vmatrix 1 & 0 & 0 & 0\\0 & 1 & 0 & 0\\
0 & 0 & -1 & 0\\0 & 0 & 0 & -1\endVmatrix
\mytag{4.6}
$$
in this frame. This frame is called an antichiral frame if 
the chirality operator $\bold H$ is given by the opposite matrix 
in this frame:
$$
\hskip -2em
H^i_j=-\Vmatrix 1 & 0 & 0 & 0\\0 & 1 & 0 & 0\\
0 & 0 & -1 & 0\\0 & 0 & 0 & -1\endVmatrix.
\mytag{4.7}
$$
\enddefinition
\mydefinition{4.3} A frame $\boldsymbol\Psi_1,\,\boldsymbol\Psi_2,
\,\boldsymbol\Psi_3,\,\boldsymbol\Psi_4$ of the Dirac bundle is 
called a self-adjoint frame if the Dirac form $\bold D$ is given
by the matrix 
$$
\hskip -2em
D_{i\bar j}=\Vmatrix 0 & 0 & 1 & 0\\0 & 0 & 0 & 1\\
1 & 0 & 0 & 0\\0 & 1 & 0 & 0\endVmatrix
\mytag{4.8}
$$
in this frame. This frame is called an anti-self-adjoint frame if 
the Dirac form $\bold D$ is given by the opposite matrix 
in this frame:
$$
\hskip -2em
D_{i\bar j}=-\Vmatrix 0 & 0 & 1 & 0\\0 & 0 & 0 & 1\\
1 & 0 & 0 & 0\\0 & 1 & 0 & 0\endVmatrix.
\mytag{4.9}
$$
\enddefinition
\mydefinition{4.4} A canonically orthonormal chiral frame 
of the Dirac bundle is a frame which is orthonormal, chiral, 
and self-adjoint simultaneously. 
\enddefinition
\mydefinition{4.5} A $P$-reverse antichiral frame of the 
Dirac bundle is a frame which is anti-orthonormal, antichiral, 
and self-adjoint simultaneously.
\enddefinition
\mydefinition{4.6} A $T$-reverse antichiral frame of the 
Dirac bundle is a frame which is orthonormal, antichiral, 
and anti-self-adjoint simultaneously.
\enddefinition
\mydefinition{4.7} A $PT$-reverse chiral frame of the 
Dirac bundle is a frame which is anti-orthonormal, chiral, 
and anti-self-adjoint simultaneously.
\enddefinition
    The definitions~\mythedefinition{4.4}, \mythedefinition{4.5}, 
\mythedefinition{4.6}, \mythedefinition{4.7} and the formulas
\mythetag{4.4}, \mythetag{4.5}, \mythetag{4.6}, \mythetag{4.7},
\mythetag{4.8}, \mythetag{4.9} describe the coordinate presentation 
of the basic field $\bold d$, $\bold H$, and $\bold D$ in special
frames whose types are listed in the diagram \mythetag{4.3}.\par
    According to the table \mythetag{4.2}, the Dirac $\gamma$-field
has one lower spacial index. Therefore its coordinate presentation 
depends not only on the choice of a spinor frame $\boldsymbol\Psi_1,
\,\boldsymbol\Psi_2,\,\boldsymbol\Psi_3,\,\boldsymbol\Psi_4$ in
the bundle $DM$, but on the choice of a tangent frame $\bold X_0,
\,\bold X_1,\,\bold X_2,\,\bold X_3$ in $TM$ too. Assume that 
$$
\boldsymbol\Psi_1,\,\boldsymbol\Psi_2,\,\boldsymbol\Psi_3,
\,\boldsymbol\Psi_4\qquad\to\qquad
\bold X_0,\,\bold X_1,\,\bold X_2,\,\bold X_3
$$
is a pair of canonically associated frames belonging to any one 
of the four types specified in the diagram \mythetag{4.3}. Then 
the components $\gamma^a_{bq}$ of the Dirac $\gamma$-field
$\boldsymbol\gamma$ are presented by the following four Dirac 
matrices in this frame pair:
$$
\xalignat 2
&\hskip -2em
\gamma^a_{b0}=\Vmatrix 0&0&1&0\\0&0&0&1\\1&0&0&0\\0&1&0&0\endVmatrix,
&&\gamma^a_{b1}=\Vmatrix 0&0&0&1\\0&0&1&0\\0&-1&0&0\\-1&0&0&0\endVmatrix,
\quad\\
\vspace{-1.5ex}
&&&\mytag{4.10}\\
\vspace{-1.5ex}
&\hskip -2em
\gamma^a_{b2}=\Vmatrix 0&0&0&-i\\0&0&i&0\\0&i&0&0\\-i&0&0&0\endVmatrix,
&&\gamma^a_{b3}=\Vmatrix 0&0&1&0\\0&0&0&-1\\-1&0&0&0\\0&1&0&0\endVmatrix.
\quad
\endxalignat
$$
The spacial index $q=0,\,1,\,2,\,3$ enumerates the matrices in 
\mythetag{4.10}. The spinor indices $a$ and $b$ determine the 
position of the component $\gamma^a_{bq}$ within one of these 
four matrices ($a$ is the row number and $b$ is the column 
number).\par
     Now let's return to the frame pairs \mythetag{3.5} and 
\mythetag{3.6}. The frames $\boldsymbol\Psi_1,\,\boldsymbol\Psi_2$
and $\tilde{\boldsymbol\Phi}_1,\,\tilde{\boldsymbol\Phi}_2$ of
the bundle $SM$ can be extended up to frames of the Dirac 
bundle $DM$. As a result we get two pairs of associated frames
$$
\align
\hskip -2em
\boldsymbol\Psi_1,\,\boldsymbol\Psi_2,
\,\boldsymbol\Psi_3,\,\boldsymbol\Psi_4\qquad&\to\qquad
\bold X_0,\,\bold X_1,\,\bold X_2,\,\bold X_3,
\mytag{4.11}\\
\vspace{1ex}
\hskip -2em
\tilde{\boldsymbol\Phi}_1,\,\tilde{\boldsymbol\Phi}_2,
\,\tilde{\boldsymbol\Phi}_3,\,\tilde{\boldsymbol\Phi}_4
\qquad&\to\qquad
\tilde\bold Y_0,\,\tilde\bold Y_1,\,\tilde\bold Y_2,
\,\tilde\bold Y_3.
\mytag{4.12}
\endalign
$$
Both frame pairs \mythetag{4.11} and \mythetag{4.12} belong 
to the first type in the diagram \mythetag{4.3}.\par
    The frames 
\mythetag{3.7} are not canonically associated. Therefore they
have no extension to a canonically associated pair. However, 
note that $\bold Y_0,\,\bold Y_1,\,\bold Y_2,\,\bold Y_3$ is
a positively polarized left orthonormal frame of $TM$ related
to the frame $\tilde\bold Y_0,\,\tilde\bold Y_1,\,\tilde\bold 
Y_2,\,\tilde\bold Y_3$ according to the formula \mythetag{1.26},
i\.\,e\. through the inversion of the space-like vectors
$\tilde\bold Y_1,\,\tilde\bold Y_2,\,\tilde\bold Y_3$. Therefore
it has an associated frame $\boldsymbol\Phi_1,
\,\boldsymbol\Phi_2,\,\boldsymbol\Phi_3,\,\boldsymbol\Phi_4$:
$$
\hskip -2em
\boldsymbol\Phi_1,\,\boldsymbol\Phi_2,
\,\boldsymbol\Phi_3,\,\boldsymbol\Phi_4
\qquad\to\qquad
\bold Y_0,\,\bold Y_1,\,\bold Y_2,\,\bold Y_3.
\mytag{4.13}
$$
The frame $\boldsymbol\Phi_1,\,\boldsymbol\Phi_2,\,\boldsymbol
\Phi_3,\,\boldsymbol\Phi_4$ in \mythetag{4.13} is produced from 
the frame $\tilde{\boldsymbol\Phi}_1,\,\tilde{\boldsymbol\Phi}_2,
\,\tilde{\boldsymbol\Phi}_3,\,\tilde{\boldsymbol\Phi}_4$
by means of the so called $P$-reversion procedure (see \mycite{3}):
$$
\xalignat 4
&\boldsymbol\Phi_1=\tilde{\boldsymbol\Phi}_3,
&&\boldsymbol\Phi_2=\tilde{\boldsymbol\Phi}_4,
&&\boldsymbol\Phi_3=\tilde{\boldsymbol\Phi}_1,
&&\boldsymbol\Phi_4=\tilde{\boldsymbol\Phi}_2.
\qquad
\mytag{4.14}
\endxalignat
$$
The frame $\boldsymbol\Phi_1,\,\boldsymbol\Phi_2,\,\boldsymbol
\Phi_3,\,\boldsymbol\Phi_4$ produced from $\tilde{\boldsymbol\Phi}_1,
\,\tilde{\boldsymbol\Phi}_2,\,\tilde{\boldsymbol\Phi}_3,\,
\tilde{\boldsymbol\Phi}_4$ by means of the formulas \mythetag{4.14}
is a $P$-reverse antichiral frame, i\.\,e\. the pair \mythetag{4.13}  
belongs to the second type of associated frame pairs in the diagram
\mythetag{4.3}.\par
     The frames $\boldsymbol\Psi_1,\,\boldsymbol\Psi_2,
\,\boldsymbol\Psi_3,\,\boldsymbol\Psi_4$ and $\tilde{\boldsymbol\Phi}_1,
\,\tilde{\boldsymbol\Phi}_2,\,\tilde{\boldsymbol\Phi}_3,
\,\tilde{\boldsymbol\Phi}_4$ in \mythetag{4.11} and \mythetag{4.12}
are related to each other by means of the formulas
$$
\xalignat 2
&\hskip -2em
\boldsymbol\Psi_i=\sum^4_{j=1}\tilde\goth T^j_i\,
\tilde{\boldsymbol\Phi}_j,
&&\tilde{\boldsymbol\Phi}_i=
\sum^4_{j=1}\tilde\goth S^j_i\,\boldsymbol\Psi_j.\quad
\mytag{4.15}
\endxalignat
$$
The formula \mythetag{4.15} is analogous to the formula 
\mythetag{2.10}. The $4\times 4$ matrices $\tilde\goth S$ 
and $\tilde\goth T$ in \mythetag{4.15} are produced from 
the matrices $\tilde\goth S$ and $\tilde\goth T$ in 
\mythetag{2.10} in some special way (see formula \thetag{2.19} 
in \mycite{3}). In our particular case they are are block-diagonal
extensions of the matrices \mythetag{2.6} and \mythetag{2.7}:
$$
\allowdisplaybreaks
\gather
\hskip -2em
\tilde\goth S=\frac{1}{|y|}
\Vmatrix
i\,y^3 & i\,y^1+y^2 & 0 & 0\\
\vspace{3ex}
i\,y^1-y^2 & -i\,y^3 & 0 & 0\\
\vspace{3ex}
0 & 0 & i\,y^3 & i\,y^1+y^2\\
\vspace{3ex}
0 & 0 & i\,y^1-y^2 & -i\,y^3
\endVmatrix,
\mytag{4.16}\\
\hskip -2em
\tilde\goth T=\frac{1}{|x|}
\Vmatrix
i\,x^3 & i\,x^1+x^2 & 0 & 0\\
\vspace{3ex}
i\,x^1-x^2 & -i\,x^3 & 0 & 0\\
\vspace{3ex}
0 & 0 & i\,x^3 & i\,x^1+x^2\\
\vspace{3ex}
0 & 0 & i\,x^1-x^2 & -i\,x^3
\endVmatrix.
\mytag{4.17}
\endgather
$$\par
     The frames $\boldsymbol\Psi_1,\,\boldsymbol\Psi_2,
\,\boldsymbol\Psi_3,\,\boldsymbol\Psi_4$ and $\boldsymbol\Phi_1,
\,\boldsymbol\Phi_2,\,\boldsymbol\Phi_3,\,\boldsymbol\Phi_4$ in 
\mythetag{4.11} and \mythetag{4.13} are also related to each 
other by means of the formulas similar to \mythetag{4.15}:
$$
\xalignat 2
&\hskip -2em
\boldsymbol\Psi_i=\sum^4_{j=1}\goth T^j_i
\,\boldsymbol\Phi_j,
&&\boldsymbol\Phi_i=\sum^4_{j=1}\goth S^j_i
\,\boldsymbol\Psi_j.\quad
\mytag{4.18}
\endxalignat
$$
The matrices $\goth S$ and $\goth T$ in \mythetag{4.18} are given
by the formulas
$$
\allowdisplaybreaks
\gather
\hskip -2em
\goth S=\frac{1}{|y|}\Vmatrix
0 & 0 & i\,y^3 & i\,y^1+y^2\\
\vspace{3ex}
0 & 0 & i\,y^1-y^2 & -i\,y^3\\
\vspace{3ex}
i\,y^3 & i\,y^1+y^2 & 0 & 0\\
\vspace{3ex}
i\,y^1-y^2 & -i\,y^3 & 0 & 0 
\endVmatrix,
\mytag{4.19}\\
\vspace{4ex}
\hskip -2em
\goth T=-\frac{1}{|x|}\Vmatrix
0 & 0 & i\,x^3 & i\,x^1+x^2\\
\vspace{3ex}
0 & 0 & i\,x^1-x^2 & -i\,x^3\\
\vspace{3ex}
i\,x^3 & i\,x^1+x^2 & 0 & 0\\
\vspace{3ex}
i\,x^1-x^2 & -i\,x^3 & 0 & 0
\endVmatrix.
\mytag{4.20}
\endgather
$$
The matrices \mythetag{4.19} and \mythetag{4.20} are similar to
\mythetag{4.17} and \mythetag{4.18}. However, unlike \mythetag{4.17} 
and \mythetag{4.18}, they are not block-diagonal. Hence, they mix 
chiral and antichiral subbundles in the expansion \mythetag{4.1}.
\head
5. Spherical coordinates.
\endhead
    The stereographic projections \mythetag{1.2} and \mythetag{1.3} 
and the local charts $x^0,\,x^1,\,x^2,\,x^3$ and $y^0,\,y^1,\,y^2,
\,y^3$ introduced through them are not very popular in physical 
literature, e\.\,g\. in \mycite{1} the spherical coordinates are
used. Therefore we consider the third chart of spherical coordinates
in $M=\Bbb R\times S^3$. The initial coordinate $\eta^0$ of these 
spherical coordinates coincides with $x^0$ and $y^0$:
$$
\pagebreak
\hskip -2em
x^0=y^0=\eta^0=\eta.
\mytag{5.1}
$$
Other three coordinates $\eta^1,\,\eta^2,\,\eta^3$ are angular 
variables:
$$
\xalignat 3
&\hskip -2em
\eta^1=\chi,
&&\eta^2=\theta,
&&\eta^3=\varphi.
\mytag{5.2}
\endxalignat
$$
They are introduced through the following formulas:
$$
\hskip -2em
\aligned
&z^1=R\,\sin\chi\,\sin\theta\,\sin\varphi,\\
&z^2=R\,\sin\chi\,\sin\theta\,\cos\varphi,\\
&z^3=R\,\sin\chi\,\cos\theta,\\
&z^4=R\,\cos\chi.
\endaligned
\mytag{5.3}
$$
The formulas relating $\eta^0,\,\eta^1,\,\eta^2,\,\eta^3$ with
$x^0,\,x^1,\,x^2,\,x^3$ and $y^0,\,y^1,\,y^2,\,y^3$ are derived 
from \mythetag{5.1} and \mythetag{5.3} by inverting \mythetag{1.2}
and \mythetag{1.3}:
$$
\xalignat 2
&\hskip -2em
\left\{\aligned
&x^0=\eta\\
&x^1=\frac{\sin\chi\,\sin\theta\,\sin\varphi}{1-\cos\chi},\\
&x^2=\frac{\sin\chi\,\sin\theta\,\cos\varphi}{1-\cos\chi},\\
&x^3=\frac{\sin\chi\,\cos\theta}{1-\cos\chi},
\endaligned\right.
&&\left\{\aligned
&y^0=\eta\\
&y^1=\frac{\sin\chi\,\sin\theta\,\sin\varphi}{1+\cos\chi},\\
&y^2=\frac{\sin\chi\,\sin\theta\,\cos\varphi}{1+\cos\chi},\\
&y^3=\frac{\sin\chi\,\cos\theta}{1+\cos\chi}.
\endaligned\right.
\mytag{5.4}
\endxalignat
$$
Differentiating the formulas \mythetag{5.4} and substituting
them into \mythetag{1.8} and \mythetag{1.9}, we derive the 
formula for the metric in the local coordinates $\eta^0,
\,\eta^1,\,\eta^2,\,\eta^3$:
$$
ds^{\kern 0.5pt 2}=R^2\,(d\eta^0)^2-
R^2\,(d\eta^1)^2-R^2\,\sin^2\!\chi\,(d\eta^2)^2
-R^2\,\sin^2\!\chi\,\sin^2\!\theta\,(d\eta^3)^2.
\quad
\mytag{5.5}
$$
The formula \mythetag{5.5} means that the holonomic coordinate frame
$$
\xalignat 4
&\hskip -2em
\frac{\partial}{\partial\eta^0},
&&\frac{\partial}{\partial\eta^1},
&&\frac{\partial}{\partial\eta^2},
&&\frac{\partial}{\partial\eta^3}
\quad
\mytag{5.6}
\endxalignat
$$
is orthogonal, but it is not an orthonormal frame. For this reason, 
instead of \mythetag{5.6}, we use the following non-holonomic 
orthonormal frame:
$$
\xalignat 2
&\hskip -2em
\bold E_0=\frac{1}{R}\,\frac{\partial}{\partial\eta^0},
&&\bold E_1=\frac{1}{R}\,\frac{\partial}{\partial\eta^1},
\quad\\
\vspace{-1ex}
&&&\mytag{5.7}\\
\vspace{-1ex}
&\hskip -2em
\bold E_2=\frac{1}{R\,\sin\chi}\,\frac{\partial}{\partial\eta^2},
&&\bold E_3=\frac{1}{R\,\sin\chi\,\sin\theta}
\,\frac{\partial}{\partial\eta^3}.\quad
\endxalignat
$$
Like in \mythetag{1.16} and \mythetag{1.18}, we have nontrivial 
commutation relationships for the frame vectors fields $\bold E_0,
\,\bold E_1,\,\bold E_2,\,\bold E_3$ defined by the formulas 
\mythetag{5.7}:
$$
[\bold E_i,\bold E_j]=\sum^3_{k=0}c^{\,k}_{ij}\,\bold E_k.
\mytag{5.8}
$$
Most of the coefficients $c^{\,k}_{ij}$ in \mythetag{5.8} are
zero. Here is the list of nonzero ones:
$$
\hskip -2em
\aligned
&c^{\,1}_{01}=-c^{\,1}_{10}=c^{\,2}_{02}=-c^{\,2}_{20}=
c^{\,3}_{03}=-c^{\,3}_{30}=-\frac{R'}{R^2},\\
\vspace{1ex}
&c^{\,2}_{12}=-c^{\,2}_{21}=c^{\,3}_{13}
=-c^{\,3}_{31}=-\frac{\cos\chi}{R\,\sin\chi},\\
\vspace{1ex}
&c^{\,3}_{23}=-c^{\,3}_{32}=-\frac{\cos\theta}
{R\,\sin\chi\,\sin\theta}.
\endaligned
\mytag{5.9}
$$
Now let's study how the frame \mythetag{5.7} is related to
the frames \mythetag{1.13} and \mythetag{1.26}:
$$
\xalignat 2
&\hskip -2em
\bold E_i=\sum^3_{j=0}\hat S^j_i\,\bold X_j,
&&\bold E_i=
\sum^3_{j=0}\Check{\Tilde S}^j_i\,\tilde\bold Y_j.\quad
\mytag{5.10}
\endxalignat
$$
By means of direct calculations, using \mythetag{5.7}, \mythetag{1.13},
\mythetag{1.14}, and \mythetag{1.26}, we derive the following formulas 
for the matrices $\hat S$ and $\Check{\Tilde S}$ in \mythetag{5.10}:
$$
\align
&\hskip -2em
\hat S=\Vmatrix 1 & 0 & 0 & 0 \\
\vspace{2ex}
0 & -\sin\varphi\,\sin\theta & \sin\varphi\,\cos\theta & \cos\varphi\\
\vspace{2ex}
0 & -\cos\varphi\,\sin\theta & \cos\varphi\,\cos\theta & -\sin\varphi\\
\vspace{2ex}
0 & -\cos\theta & -\sin\theta & 0
\endVmatrix,
\mytag{5.11}\\
\vspace{3ex}
&\hskip -2em
\Check{\Tilde S}=\Vmatrix 1 & 0 & 0 & 0 \\
\vspace{2ex}
0 & -\sin\varphi\,\sin\theta & -\sin\varphi\,\cos\theta & -\cos\varphi\\
\vspace{2ex}
0 & -\cos\varphi\,\sin\theta & -\cos\varphi\,\cos\theta & \sin\varphi\\
\vspace{2ex}
0 & -\cos\theta & \sin\theta & 0
\endVmatrix.
\mytag{5.12}
\endalign
$$
As appears, both matrices \mythetag{5.11} and \mythetag{5.12} belong 
to the special orthochronous Lorentz group $\MatGrSO^+(1,3,\Bbb R)$. 
For this reason they are related with two matrices $\hat\goth S$ and
$\Check{\Tilde{\goth S}}$ belonging to the group $\MatGrSL(2,\Bbb C)$
through the group homomorphism \mythetag{2.1}:
$$
\xalignat 2
&\hskip -2em
\hat S=\phi(\hat\goth S),
&&\Check{\Tilde S}=\phi(\Check{\Tilde{\goth S}}).
\mytag{5.13}
\endxalignat
$$
These two matrices $\hat\goth S$ and $\Check{\Tilde{\goth S}}$ 
satisfying \mythetag{5.13} can be found in explicit form:
$$
\align
&\hskip -2em
\hat\goth S=\frac{1}{\sqrt{2}}
\Vmatrix \exp\Bigl(\dfrac{i\,\varphi+i\,\theta}{2}\Bigr)
& -\exp\Bigl(\dfrac{i\,\varphi-i\,\theta}{2}\Bigr)\\
\vspace{4ex}
\exp\Bigl(\dfrac{i\,\theta-i\,\varphi}{2}\Bigr) &
\exp\Bigl(\dfrac{-i\,\theta-i\,\varphi}{2}\Bigr)
\endVmatrix\quad
\mytag{5.14}\\
\vspace{4ex}
&\hskip -2em
\Check{\Tilde{\goth S}}=
\frac{i}{\sqrt{2}}
\Vmatrix -\exp\Bigl(\dfrac{i\,\varphi-i\,\theta}{2}\Bigr)
& \exp\Bigl(\dfrac{i\,\varphi+i\,\theta}{2}\Bigr)\\
\vspace{4ex}
\exp\Bigl(\dfrac{-i\,\theta-i\,\varphi}{2}\Bigr) &
\exp\Bigl(\dfrac{i\,\theta-i\,\varphi}{2}\Bigr)
\endVmatrix\quad
\mytag{5.15}
\endalign
$$
If we express the matrix \mythetag{2.6} through the spherical coordinates
\mythetag{5.2}, then we get
$$
\tilde\goth S=
\Vmatrix i\,\cos\theta & \sin\theta\,\exp(i\varphi)\\
\vspace{4ex}
-\sin\theta\,\exp(-i\varphi) & -i\,\cos\theta
\endVmatrix.\quad
\mytag{5.16}
$$
The matrices \mythetag{5.14}, \mythetag{5.15}, and \mythetag{5.16} are
related to each other as follows:
$$
\hskip -2em
\hat\goth S=\tilde\goth S\,\Check{\Tilde{\goth S}}.
\mytag{5.17}
$$
Applying the group homomorphism \mythetag{2.1} to \mythetag{5.17}, we
get:
$$
\hskip -2em
\hat S=\tilde S\,\Check{\Tilde S}.
\mytag{5.18}
$$
The relationship \mythetag{5.18} can be verified directly using
the formulas \mythetag{5.11}, \mythetag{5.12}, and the the
formula for $\hat S$ in section 1. Now we define a spinor frame
$\boldsymbol\Xi_1,\,\boldsymbol\Xi_2$ of the bundle $SM$
associated with the tangent frame $\bold E_0,\,\bold E_1,
\,\bold E_2,\,\bold E_3$. We do it by setting
$$
\xalignat 2
&\hskip -2em
\boldsymbol\Xi_i=\sum^2_{j=1}\hat\goth S^j_i\,\boldsymbol\Psi_j,
&&\boldsymbol\Xi_i=
\sum^2_{j=1}\Check{\Tilde{\goth S}}^j_i\,\tilde{\boldsymbol\Phi}_j.
\quad
\mytag{5.19}
\endxalignat
$$
Due to \mythetag{5.17}, \mythetag{5.10}, and \mythetag{2.10} the
formulas \mythetag{5.19} are consistent with each other. The frame
$\boldsymbol\Xi_1,\,\boldsymbol\Xi_2$ determined by any one of these
two formulas is an orthonormal frame of the bundle of Weyl spinors.
So we have the frame association 
$$
\hskip -2em
\boldsymbol\Xi_1,\,\boldsymbol\Xi_2\qquad\to\qquad
\bold E_0,\,\bold E_1,\,\bold E_2,\,\bold E_3
\mytag{5.20}
$$
analogous to \mythetag{3.5} and \mythetag{3.6}. This frame association
\mythetag{5.20} is a special case of the general scheme presented by
the diagram \mythetag{3.3}.\par
    The orthonormal frame $\boldsymbol\Xi_1,\,\boldsymbol\Xi_2$ of the
bundle of Weyl spinors $SM$ has a unique extension $\boldsymbol\Xi_1,
\,\boldsymbol\Xi_2,\,\boldsymbol\Xi_3,\,\boldsymbol\Xi_4$ to the 
Dirac bundle $DM$. This extension is a canonically orthonormal chiral 
frame associated with the frame $\bold E_0,\,\bold E_1,\,\bold E_2,\,
\bold E_3$:
$$
\hskip -2em
\boldsymbol\Xi_1,\,\boldsymbol\Xi_2,\,\boldsymbol\Xi_3,
\,\boldsymbol\Xi_4\qquad\to\qquad
\bold E_0,\,\bold E_1,\,\bold E_2,\,\bold E_3.
\mytag{5.21}
$$
The frame association \mythetag{5.21} is analogous to \mythetag{4.11}
and \mythetag{4.12}. Note that the frame $\bold E_0,\,\bold E_1,
\,\bold E_2,\,\bold E_3$ in \mythetag{5.20} and \mythetag{5.21} is
a positively polarized right orthonormal frame. This fact follows from
$\hat S\in\MatGrSO^+(1,3,\Bbb R)$ and from \mythetag{5.10}. Therefore
the frame association \mythetag{5.21} is a special case for the scheme 
presented by the first line in the diagram \mythetag{4.3}.\par
    The extended frame $\boldsymbol\Xi_1,\,\boldsymbol\Xi_2,
\,\boldsymbol\Xi_3,\,\boldsymbol\Xi_4$ is related to the previously
defined frames $\boldsymbol\Psi_1,\,\boldsymbol\Psi_2,\,\boldsymbol
\Psi_3,\,\boldsymbol\Psi_4$ and $\tilde{\boldsymbol\Phi}_1,
\,\tilde{\boldsymbol\Phi}_2,\,\tilde{\boldsymbol\Phi}_3,
\,\tilde{\boldsymbol\Phi}_4$ as follows
$$
\xalignat 2
&\hskip -2em
\boldsymbol\Xi_i=\sum^4_{j=1}\hat\goth S^j_i\,\boldsymbol\Psi_j,
&&\boldsymbol\Xi_i=
\sum^4_{j=1}\Check{\Tilde{\goth S}}^j_i\,\tilde{\boldsymbol\Phi}_j.
\quad
\mytag{5.22}
\endxalignat
$$
The matrices $\hat\goth S$ and $\Check{\Tilde{\goth S}}$ 
in \mythetag{5.22} are four-dimensional extensions of the matrices 
\mythetag{5.14} and \mythetag{5.14}. They are block-diagonal 
matrices constructed with the use of the two-dimensional matrices
$\hat\goth S$ and $\Check{\Tilde{\goth S}}$. Here are the explicit
formulas for them:
$$
\allowdisplaybreaks
\align
&\hskip -2em
\hat\goth S=\frac{1}{\sqrt{2}}
\Vmatrix e^{\raise 4pt\hbox{$\frac{i\varphi+i\theta}{2}$}}
& -e^{\raise 4pt\hbox{$\frac{i\,\varphi-i\,\theta}{2}$}} & 0 & 0\\
\vspace{4ex}
e^{\raise 4pt\hbox{$\frac{i\,\theta-i\,\varphi}{2}$}}  &
e^{\raise 4pt\hbox{$\frac{-i\,\theta-i\,\varphi}{2}$}} & 0 & 0\\
\vspace{4ex}
0 & 0 &e^{\raise 4pt\hbox{$\frac{i\,\varphi+i\,\theta}{2}$}}
& -e^{\raise 4pt\hbox{$\frac{i\,\varphi-i\,\theta}{2}$}}\\
\vspace{4ex}
0 & 0 &e^{\raise 4pt\hbox{$\frac{i\,\theta-i\,\varphi}{2}$}} & 
e^{\raise 4pt\hbox{$\frac{-i\,\theta-i\,\varphi}{2}$}}
\endVmatrix,\quad
\mytag{5.23}\\
\vspace{4ex}
&\hskip -2em
\Check{\Tilde{\goth S}}=
\frac{i}{\sqrt{2}}
\Vmatrix -e^{\raise 4pt\hbox{$\frac{i\,\varphi-i\,\theta}{2}$}}
& e^{\raise 4pt\hbox{$\frac{i\,\varphi+i\,\theta}{2}$}} & 0 & 0\\
\vspace{4ex}
e^{\raise 4pt\hbox{$\frac{-i\,\theta-i\,\varphi}{2}$}} &
e^{\raise 4pt\hbox{$\frac{i\,\theta-i\,\varphi}{2}$}} & 0 & 0\\
\vspace{4ex}
0 & 0 &-e^{\raise 4pt\hbox{$\frac{i\,\varphi-i\,\theta}{2}$}}
& e^{\raise 4pt\hbox{$\frac{i\,\varphi+i\,\theta}{2}$}}\\
\vspace{4ex}
0 & 0 &e^{\raise 4pt\hbox{$\frac{-i\,\theta-i\,\varphi}{2}$}} &
e^{\raise 4pt\hbox{$\frac{i\,\theta-i\,\varphi}{2}$}}
\endVmatrix.\quad
\mytag{5.24}\\
\endalign
$$
The matrix \mythetag{5.16} has an analogous four-dimensional extension.
Here is the explicit formula for such an extension:
$$
\hskip -2em
\tilde\goth S=
\Vmatrix i\,\cos\theta & \sin\theta\,e^{\raise 2pt
\hbox{$\ssize i\varphi$}} & 0 & 0\\
\vspace{4ex}
-\sin\theta\,e^{\raise 2pt\hbox{$\ssize -i\varphi$}} 
& -i\,\cos\theta & 0 & 0\\
\vspace{4ex}
0 & 0 & i\,\cos\theta & \sin\theta\,e^{\raise 2pt
\hbox{$\ssize i\varphi$}}\\
\vspace{4ex}
0 & 0 & -\sin\theta\,e^{\raise 2pt\hbox{$\ssize -i\varphi$}} 
& -i\,\cos\theta
\endVmatrix.\quad
\mytag{5.25}
$$
The matrix \mythetag{5.25} is also a block-diagonal matrix. It is
used in the following transition formula relating the frames
$\boldsymbol\Psi_1,\,\boldsymbol\Psi_2,\,\boldsymbol
\Psi_3,\,\boldsymbol\Psi_4$ and $\tilde{\boldsymbol\Phi}_1,
\,\tilde{\boldsymbol\Phi}_2,\,\tilde{\boldsymbol\Phi}_3,
\,\tilde{\boldsymbol\Phi}_4$:
$$
\tilde{\boldsymbol\Phi}_i=
\sum^4_{j=1}\tilde\goth S^j_i\,\boldsymbol\Psi_j.
\mytag{5.26}
$$
The formula \mythetag{5.26} is an extension for the second formula
\mythetag{2.10}. It coincides with the second formula \mythetag{4.15}. 
As for the matrices \mythetag{5.23}, \mythetag{5.24}, and \mythetag{5.25},
being extensions of the matrices \mythetag{5.14}, \mythetag{5.15}, 
and \mythetag{5.16}, they satisfy the relationship \mythetag{5.17}. 
The matrix \mythetag{5.25} can be produced from \mythetag{4.16} 
by passing to the angular variables \mythetag{5.2}.\par
    Let's recall that the frame $\bold Y_0,\,\bold Y_1,\,\bold Y_2,
\,\bold Y_3$ has no canonically associated spinor frame in $SM$ (since
it is left oriented), but has an associated frame $\boldsymbol\Phi_1,
\,\boldsymbol\Phi_2,\,\boldsymbol\Phi_3,\,\boldsymbol\Phi_4$ in $DM$ 
(see \mythetag{4.13} above). For this reason we can write
$$
\hskip -2em
\boldsymbol\Xi_i=
\sum^4_{j=1}\check{\goth S}^j_i\,\boldsymbol\Phi_j.
\mytag{5.27}
$$
The matrix $\check\goth S$ in \mythetag{5.27} is given by the following
formula:
$$
\hskip -2em
\check\goth S=
\frac{i}{\sqrt{2}}
\Vmatrix 0 & 0 & -e^{\raise 4pt\hbox{$\frac{i\,\varphi-i\,\theta}{2}$}}
& e^{\raise 4pt\hbox{$\frac{i\,\varphi+i\,\theta}{2}$}}\\
\vspace{4ex}
0 & 0 & e^{\raise 4pt\hbox{$\frac{-i\,\theta-i\,\varphi}{2}$}} &
e^{\raise 4pt\hbox{$\frac{i\,\theta-i\,\varphi}{2}$}}\\
\vspace{4ex}
-e^{\raise 4pt\hbox{$\frac{i\,\varphi-i\,\theta}{2}$}}
& e^{\raise 4pt\hbox{$\frac{i\,\varphi+i\,\theta}{2}$}} &0 & 0\\
\vspace{4ex}
e^{\raise 4pt\hbox{$\frac{-i\,\theta-i\,\varphi}{2}$}} &
e^{\raise 4pt\hbox{$\frac{i\,\theta-i\,\varphi}{2}$}} & 0 & 0
\endVmatrix.\quad
\mytag{5.28}
$$
As we see, the matrix \mythetag{5.28} is not block-diagonal. It is 
because the frame pair \mythetag{4.13} corresponds to the second
line in the diagram \mythetag{4.3}.
\head
6. Metric connection and its spinor components.
\endhead
     In previous sections we considered three local chart in the
spherical universe $M=\Bbb R\times S^3$ --- two charts with 
stereographic coordinates and one chart with spherical coordinates.
Considering Dirac spinors, we have equipped these local charts with
three associated frame pairs (see \mythetag{4.11}, \mythetag{4.13},
and \mythetag{5.21} above). Our next goal is to calculate the
components of the metric connection $(\Gamma,\Alpha,\bar{\Alpha})$
in each of these three frame pairs.\par
     Spinor extension of the metric connection has three groups of
components: spacial components $\Gamma^k_{ij}$, spinor components 
$\Alpha^b_{i\kern 1pt a}$, and conjugate spinor components 
$\bar{\Alpha}\vphantom{\Alpha}^{\bar a}_{i\kern 1pt \bar b}$. 
Assume for a while that we have some arbitrary frame pair
$$
\hskip -2em
\boldsymbol\Psi_1,\,\boldsymbol\Psi_2,\,\boldsymbol\Psi_3,
\,\boldsymbol\Psi_4\qquad-\qquad
\bold X_0,\,\bold X_1,\,\bold X_2,\,\bold X_3,
\mytag{6.1}
$$
where $\boldsymbol\Psi_1,\,\boldsymbol\Psi_2,\,\boldsymbol\Psi_3,
\,\boldsymbol\Psi_4$ is a spinor frame in $DM$ and $\bold X_0,\,
\bold X_1,\,\bold X_2,\,\bold X_3$ is a spacial frame in $TM$.
They can be either associated or non-associated frames. In any case
the conjugate spinor components $\bar{\Alpha}\vphantom{\Alpha}^{\bar
a}_{i\kern 1pt \bar b}$ are expressed through spinor components in
an elementary way through complex conjugation: 
$$
\hskip -2em
\bar{\Alpha}\vphantom{\Alpha}^{\bar a}_{i\kern 1pt \bar b}
=\overline{\Alpha^{\bar a}_{i\kern 1pt \bar b}}.
\mytag{6.2}
$$
Due to the formula \mythetag{6.2} it is sufficient to know 
$\Gamma^k_{ij}$ and $\Alpha^b_{i\kern 1pt a}$ for to describe 
the metric connection completely. The spacial components 
$\Gamma^k_{ij}$ correspond to the well-known Levi-Civita
connection in $TM$. They are given by the formula
$$
\hskip -2em
\gathered
\Gamma^k_{ij}=\sum^3_{r=0}\frac{g^{\kern 0.5pt kr}}{2}
\left(L_{\bold X_i}\!(g_{rj})
+L_{\bold X_j}\!(g_{i\kern 0.5pt r})
-L_{\bold X_r}\!(g_{ij})\right)+\\
+\,\frac{c^{\,k}_{ij}}{2}
-\sum^3_{r=0}\sum^3_{s=0}\frac{c^{\,s}_{i\kern 0.5pt r}}{2}\,g^{kr}
\,g_{sj}-\sum^3_{r=0}\sum^3_{s=0}\frac{c^{\,s}_{j\kern 0.5ptr}}{2}
\,g^{kr}\,g_{s\kern 0.5pt i}
\endgathered
\mytag{6.3}
$$
(see \mycite{4}). The Levi-Civita connection is a torsion-free
connection. Nevertheless, its components $\Gamma^k_{ij}$ in 
\mythetag{6.3} are not symmetric: $\Gamma^k_{ij}-\Gamma^k_{j\kern 
0.5pt i}=c^{\,k}_{ij}$. It is because in general case the tangent 
frame $\bold X_0,\,\bold X_1,\,\bold X_2,\,\bold X_3$ in the pair 
\mythetag{6.1} is not commutative (see \mythetag{1.16} above).\par
    Now let's proceed to the spinor components $\Alpha^b_{i\kern 1pt a}$
of our metric connection. These components are given
by the formula
$$
\gathered
\Alpha^a_{ib}
=\sum^4_{\alpha=1}\sum^4_{\beta=1}\frac{
L_{\bold X_i}(d_{\alpha\beta})
\,d^{\kern 0.5pt \beta\kern 0.5pt\alpha}}{8}\,\delta^a_b
-\sum^4_{\alpha=1}\sum^4_{\beta=1}\sum^4_{d=1}\frac{
L_{\bold X_i}(d_{\alpha\beta})
\,d^{\kern 0.5pt \beta\kern 0.5pt d}\,H^\alpha_d}{8}\,H^a_b\,-\\
-\sum^4_{c=1}\sum^4_{d=1}
\sum^4_{r=1}\frac{d_{\kern 0.5pt bc}
\,L_{\bold X_i}(H^c_d)
\,H^d_r\,d^{\kern 0.5pt ra}}{4}\,+
\sum^3_{m=0}\sum^3_{n=0}\sum^4_{\alpha=1}
\frac{L_{\bold X_i}(\gamma^{\,\alpha}_{b\kern 0.5pt m}
\,g^{mn})}{4}\,\times\\
\times\,\gamma^{\,a}_{\alpha\kern 0.5pt n}
+\sum^3_{m=0}\sum^3_{n=0}\sum^4_{\alpha=1}\sum^3_{s=0}
\frac{\gamma^{\,\alpha}_{b\kern 0.5pt m}\,\Gamma^n_{is}\,g^{ms}
\,\gamma^{\,a}_{\alpha\kern 0.5pt n}}{4}
\endgathered\qquad
\mytag{6.4}
$$
(see \mycite{4} for more details). This is a general formula applicable
to an arbitrary frame pair \mythetag{6.1}. In our particular case all
of our three frame pairs, which we study in this section, are special
ones. They are canonically associated frame pairs described by the first
and second lines in the diagram \mythetag{4.3}. In each such frame pair
the components of the basic fields $\bold g$, $\bold d$, $\bold H$, 
$\bold D$, and $\boldsymbol\gamma$ are constants (see formulas 
\mythetag{1.15}, \mythetag{4.4}, \mythetag{4.5}, \mythetag{4.6}, 
\mythetag{4.7}, \mythetag{4.8}, \mythetag{4.9}, and \mythetag{4.10}).
Therefore, their derivatives $L_{\bold X_i}$ are identically zero. 
As a result the formulas \mythetag{6.3} and \mythetag{6.4} are
reduced to
$$
\align
&\hskip -2em
\Gamma^k_{ij}=\frac{c^{\,k}_{ij}}{2}
-\sum^3_{r=0}\sum^3_{s=0}\frac{c^{\,s}_{i\kern 0.5pt r}}{2}\,g^{kr}
\,g_{sj}-\sum^3_{r=0}\sum^3_{s=0}\frac{c^{\,s}_{j\kern 0.5ptr}}{2}
\,g^{kr}\,g_{s\kern 0.5pt i},
\mytag{6.5}\\
\vspace{2ex}
&\hskip -2em
\Alpha^a_{ib}=\sum^3_{m=0}\sum^3_{n=0}\sum^4_{\alpha=1}\sum^3_{s=0}
\frac{\gamma^{\,\alpha}_{b\kern 0.5pt m}\,\Gamma^n_{is}\,g^{ms}
\,\gamma^{\,a}_{\alpha\kern 0.5pt n}}{4}.
\mytag{6.6}
\endalign
$$
The metric connection with the components \mythetag{6.5} and
\mythetag{6.6} produces two curvature tensors. Here are the
formulas for their components:
$$
\align
&\hskip -4em
R^p_{qij}=L_{\bold X_i}(\Gamma^p_{\!j\,q})
-L_{\bold X_j}(\Gamma^p_{\!i\,q})
+\sum^3_{h=0}\left(\Gamma^p_{\!i\,h}\,\Gamma^h_{\!j\,q}
-\Gamma^p_{\!j\,h}\,\Gamma^h_{\!i\,q}\right)
-\sum^3_{k=0}c^{\,k}_{ij}\,\Gamma^p_{kq},
\mytag{6.7}
\hskip -1em
\\
&\hskip -4em
\goth R^p_{qij}=L_{\bold X_i}(\Alpha^p_{j\,q})
-L_{\bold X_j}(\Alpha^p_{i\,q})
+\sum^4_{h=1}\left(\Alpha^p_{i\,h}\,\Alpha^{\!h}_{j\,q}
-\Alpha^p_{j\,h}\,\Alpha^{\!h}_{i\,q}\right)
-\sum^3_{k=0}c^{\,k}_{ij}\,\Alpha^p_{kq}.
\hskip -1em
\mytag{6.8}
\endalign
$$
The formula \mythetag{6.7} yields the components of the well-known
Riemannian curvature tensor, while \mythetag{6.8} are the components
of its spinor extension. As appears, the tensors $\bold R$ and 
$\eufb R$ are related to each other as follows:
$$
\hskip -2em
\goth R^p_{qij}=\frac{1}{4}\sum^3_{m=0}\sum^3_{n=0}\sum^3_{r=0}
\sum^4_{\alpha=1}R^r_{mij}\,\gamma^\alpha_{qn}\ g^{mn}
\,\gamma^p_{\alpha r}
\mytag{6.9}
$$
(see proof of the formula \mythetag{6.9} in \mycite{4}). Unlike
\mythetag{6.5} and \mythetag{6.6}, the formulas \mythetag{6.7},
\mythetag{6.8}, and \mythetag{6.9} are applicable to an arbitrary
frame pair \mythetag{6.1}.\par
     The next step now is to apply the above formulas to our
three local charts and three frame pairs and get some formulas
specific to the homogeneous and isotropic spherical universe 
$M=\Bbb R\times S^3$.\par
     1. {\bf Stereographic coordinates} in projection from the 
{\bf North Pole}. We denote these coordinates with $x^0,\,x^1,
\,x^2,\,x^3$. Their domain is the whole sphere $S^3$
except for the North Pole itself. The frames \mythetag{4.11} 
are defined and smooth in this domain. This frame pair is linked 
to the coordinates $x^0,\,x^1,\,x^2,\,x^3$. The commutation 
coefficients for the frame $\bold X_0,\,\bold X_1,\,\bold X_2,
\,\bold X_3$ in this pair are given by the formulas
\mythetag{1.17}. Now we substitute these coefficients into
the formula \mythetag{6.5}. As a result we get the following 
complete list of nonzero $\Gamma$-components of the metric
connection in the non-holonomic frame $\bold X_0,\,\bold X_1,
\,\bold X_2,\,\bold X_3$:
$$
\xalignat 3
&\hskip -2em
\Gamma^0_{11}=\frac{R'}{R^2},
&&\Gamma^0_{22}=\frac{R'}{R^2},
&&\Gamma^0_{33}=\frac{R'}{R^2},\\
\vspace{1ex}
&\hskip -2em
\Gamma^1_{10}=\frac{R'}{R^2},
&&\Gamma^2_{20}=\frac{R'}{R^2},
&&\Gamma^3_{30}=\frac{R'}{R^2},\\
\vspace{1ex}
&\hskip -2em
\Gamma^1_{22}=\frac{(x^1)}{R^2},
&&\Gamma^2_{33}=\frac{(x^2)}{R^2},
&&\Gamma^3_{11}=\frac{(x^3)}{R^2},\\
\vspace{-1ex}
&&&\mytag{6.10}\\
\vspace{-1ex}
&\hskip -2em
\Gamma^2_{11}=\frac{(x^2)}{R^2},
&&\Gamma^3_{22}=\frac{(x^3)}{R^2},
&&\Gamma^1_{33}=\frac{(x^1)}{R^2},\\
\vspace{1ex}
&\hskip -2em
\Gamma^1_{12}=-\frac{(x^2)}{R^2},
&&\Gamma^2_{23}=-\frac{(x^3)}{R^2},
&&\Gamma^3_{31}=-\frac{(x^1)}{R^2},\\
\vspace{1ex}
&\hskip -2em
\Gamma^2_{21}=-\frac{(x^1)}{R^2},
&&\Gamma^3_{32}=-\frac{(x^2)}{R^2},
&&\Gamma^1_{13}=-\frac{(x^2)}{R^2}.
\endxalignat
$$
Here $R'$ is the derivative of the function \mythetag{1.10}.
Substituting \mythetag{6.10} into \mythetag{6.6}, we derive 
the explicit formulas for $\Alpha$-components of the metric
connection:
$$
\allowdisplaybreaks
\xalignat 3
&\hskip -2em
\Alpha^1_{11}=-\frac{i\,(x^2)}{2\,R},
&&\Alpha^1_{21}=\frac{i\,(x^1)}{2\,R},
&&\Alpha^1_{31}=\frac{R'}{2\,R^2},
\quad\\
\vspace{1ex}
&\hskip -2em
\Alpha^2_{12}=\frac{i\,(x^2)}{2\,R},
&&\Alpha^2_{22}=-\frac{i\,(x^1)}{2\,R},
&&\Alpha^2_{32}=-\frac{R'}{2\,R^2},
\quad\\
\vspace{-1ex}
&&&\mytag{6.11}\\
\vspace{-1ex}
&\hskip -2em
\Alpha^3_{13}=-\frac{i\,(x^2)}{2\,R},
&&\Alpha^3_{23}=\frac{i\,(x^1)}{2\,R},
&&\Alpha^3_{33}=-\frac{R'}{2\,R^2},
\quad\\
\vspace{1ex}
&\hskip -2em
\Alpha^4_{14}=\frac{i\,(x^2)}{2\,R},
&&\Alpha^4_{24}=-\frac{i\,(x^1)}{2\,R},
&&\Alpha^4_{34}=\frac{R'}{2\,R^2},
\quad\\
\endxalignat
$$
$$
\allowdisplaybreaks
\xalignat 2
&\hskip -2em
\Alpha^1_{12}=\frac{R'}{2\,R^2}+\frac{(x^3)}{2\,R},
&&\Alpha^2_{11}=\frac{R'}{R^2}-\frac{(x^3)}{2\,R},
\quad\\
\vspace{1ex}
&\hskip -2em
\Alpha^3_{14}=-\frac{R'}{2\,R^2}+\frac{(x^3)}{2\,R},
&&\Alpha^4_{13}=-\frac{R'}{2\,R^2}-\frac{(x^3)}{2\,R},
\quad\\
\vspace{1ex}
&\hskip -2em
\Alpha^1_{22}=-\frac{i\,R'}{2\,R^2}-\frac{i\,(x^3)}{2\,R},
&&\Alpha^2_{21}=\frac{i\,R'}{2\,R^2}-\frac{i\,(x^3)}{2\,R},
\quad\\
\vspace{-1ex}
&&&\mytag{6.12}
\quad\\
\vspace{-1ex}
&\hskip -2em
\Alpha^3_{24}=\frac{i\,R'}{2\,R^2}-\frac{i\,(x^3)}{2\,R},
&&\Alpha^4_{23}=-\frac{i\,R'}{2\,R^2}-\frac{i\,(x^3)}{2\,R},
\quad\\
\vspace{1ex}
&\hskip -2em
\Alpha^1_{32}=-\frac{(x^1)}{2\,R}+\frac{i\,(x^2)}{2\,R},
&&\Alpha^2_{31}=\frac{(x^1)}{2\,R}+\frac{i\,(x^2)}{2\,R},
\quad\\
\vspace{1ex}
&\hskip -2em
\Alpha^3_{34}=-\frac{(x^1)}{2\,R}+\frac{i\,(x^2)}{2\,R},
&&\Alpha^4_{33}=\frac{(x^1)}{2\,R}+\frac{i\,(x^2)}{2\,R}.
\quad
\endxalignat
$$
Substituting \mythetag{6.10} into \mythetag{6.7}, we find the
components of the Riemannian curvature tensor $\bold R$. Its 
nonzero components are listed here:
$$
\xalignat 2
&\hskip -1ex
R^0_{101}=-R^0_{110}=-\frac{(R')^2}{R^4}+\frac{R''}{R^3},
&&R^1_{001}=-R^1_{010}=-\frac{(R')^2}{R^4}+\frac{R''}{R^3},
\qquad\quad\\
\vspace{1ex}
&\hskip -1ex
R^0_{202}=-R^0_{220}=-\frac{(R')^2}{R^4}+\frac{R''}{R^3},
&&R^2_{002}=-R^2_{020}=-\frac{(R')^2}{R^4}+\frac{R''}{R^3},
\qquad\quad\\
\vspace{1ex}
&\hskip -1ex
R^0_{303}=-R^0_{330}=-\frac{(R')^2}{R^4}+\frac{R''}{R^3},
&&R^3_{003}=-R^3_{030}=-\frac{(R')^2}{R^4}+\frac{R''}{R^3},
\qquad\quad\\
\vspace{-1ex}
&&&\qquad\quad
\mytag{6.13}\\
\vspace{-1ex}
&\hskip -1ex
R^1_{212}=-R^1_{221}=\frac{1}{R^2}+\frac{(R')^2}{R^4},
&&R^2_{112}=-R^2_{121}=-\frac{1}{R^2}-\frac{(R')^2}{R^4},
\qquad\quad\\
\vspace{1ex}
&\hskip -1ex
R^2_{323}=-R^2_{333}=\frac{1}{R^2}+\frac{(R')^2}{R^4},
&&R^3_{223}=-R^3_{232}=-\frac{1}{R^2}-\frac{(R')^2}{R^4},
\qquad\quad\\
\vspace{1ex}
&\hskip -1ex
R^3_{131}=-R^3_{113}=\frac{1}{R^2}+\frac{(R')^2}{R^4},
&&R^1_{331}=-R^1_{313}=-\frac{1}{R^2}-\frac{(R')^2}{R^4}.
\qquad\quad
\endxalignat
$$
The next step is to substitute \mythetag{6.11} and \mythetag{6.12}
into \mythetag{6.8}. As a result we derive the explicit formulas 
for the components of the spinor curvature tensor $\eufb R$. Below 
is the list of all its nonzero components:
$$
\allowdisplaybreaks
\xalignat 2
&\hskip -1.5em
\goth R^1_{201}=-\goth R^1_{210}
=-\frac{(R')^2}{2\,R^4}+\frac{R''}{2\,R^3},
&&\goth R^2_{101}=-\goth R^2_{110}
=-\frac{(R')^2}{2\,R^4}+\frac{R''}{2\,R^3},
\qquad\\
\vspace{1ex}
&\hskip -1.5em
\goth R^3_{401}=-\goth R^3_{410}
=\frac{(R')^2}{2\,R^4}-\frac{R''}{2\,R^3},
&&\goth R^4_{301}=-\goth R^4_{310}
=-\frac{(R')^2}{2\,R^4}+\frac{R''}{2\,R^3},
\qquad\\
\vspace{1ex}
&\hskip -1.5em
\goth R^1_{202}=-\goth R^1_{220}
=\frac{i\,(R')^2}{2\,R^4}-\frac{i\,R''}{2\,R^3},
&&\goth R^2_{102}=-\goth R^2_{120}
=-\frac{i\,(R')^2}{2\,R^4}+\frac{i\,R''}{2\,R^3},
\qquad
\mytag{6.14}\\
\vspace{1ex}
&\hskip -1.5em
\goth R^3_{402}=-\goth R^3_{420}
=-\frac{i\,(R')^2}{2\,R^4}+\frac{i\,R''}{2\,R^3},
&&\goth R^4_{302}=-\goth R^4_{320}
=\frac{i\,(R')^2}{2\,R^4}-\frac{i\,R''}{2\,R^3},
\qquad\\
\vspace{1ex}
&\hskip -1.5em
\goth R^1_{103}=-\goth R^1_{130}
=-\frac{(R')^2}{2\,R^4}+\frac{R''}{2\,R^3},
&&\goth R^2_{203}=-\goth R^2_{230}
=\frac{(R')^2}{2\,R^4}-\frac{R''}{2\,R^3},
\qquad\\
\vspace{1ex}
&\hskip -1.5em
\goth R^3_{303}=-\goth R^3_{330}
=\frac{(R')^2}{2\,R^4}-\frac{R''}{2\,R^3},
&&\goth R^4_{403}=-\goth R^4_{430}
=-\frac{(R')^2}{2\,R^4}+\frac{R''}{2\,R^3},
\qquad\\
\vspace{1ex}
&\hskip -1.5em
\goth R^1_{112}=-\goth R^1_{121}
=\frac{i}{2\,R^2}+\frac{i\,(R')^2}{2\,R^4},
&&\goth R^2_{212}=-\goth R^2_{221}
=-\frac{i}{2\,R^2}-\frac{i\,(R')^2}{2\,R^4},
\qquad\\
\vspace{1ex}
&\hskip -1.5em
\goth R^3_{312}=-\goth R^3_{321}
=\frac{i}{2\,R^2}+\frac{i\,(R')^2}{2\,R^4},
&&\goth R^4_{412}=-\goth R^4_{421}
=-\frac{i}{2\,R^2}-\frac{i\,(R')^2}{2\,R^4},
\qquad\\
\vspace{1ex}
&\hskip -1.5em
\goth R^1_{223}=-\goth R^1_{232}
=\frac{i}{2\,R^2}+\frac{i\,(R')^2}{2\,R^4},
&&\goth R^2_{123}=-\goth R^2_{132}
=\frac{i}{2\,R^2}+\frac{i\,(R')^2}{2\,R^4},
\qquad\mytag{6.15}\\
\vspace{1ex}
&\hskip -1.5em
\goth R^3_{423}=-\goth R^3_{432}
=\frac{i}{2\,R^2}+\frac{i\,(R')^2}{2\,R^4},
&&\goth R^4_{323}=-\goth R^4_{332}
=\frac{i}{2\,R^2}+\frac{i\,(R')^2}{2\,R^4},
\qquad\\
\vspace{1ex}
&\hskip -1.5em
\goth R^1_{231}=-\goth R^1_{213}
=\frac{1}{2\,R^2}+\frac{(R')^2}{2\,R^4},
&&\goth R^2_{131}=-\goth R^2_{131}
=-\frac{1}{2\,R^2}-\frac{(R')^2}{2\,R^4},
\qquad\\
\vspace{1ex}
&\hskip -1.5em
\goth R^3_{431}=-\goth R^3_{413}
=\frac{1}{2\,R^2}+\frac{(R')^2}{2\,R^4},
&&\goth R^4_{331}=-\goth R^4_{331}
=-\frac{1}{2\,R^2}-\frac{(R')^2}{2\,R^4}.
\qquad
\endxalignat
$$
One can verify the relationship \mythetag{6.9} by direct
calculations as a test for consistency of \mythetag{6.13}, 
\mythetag{6.14}, and \mythetag{6.15}. The Ricci tensor is
an important object in general relativity. It is used in
Einstein's equation (see \mycite{1} or \mycite{5}). Using
\mythetag{6.13}, we can calculate the components of the
Ricci tensor in the frame $\bold X_0,\,\bold X_1,
\,\bold X_2,\,\bold X_3$. As appears, the matrix of the 
Ricci tensor is diagonal. Here are its diagonal elements 
$$
\xalignat 2
&\hskip -2em
R_{00}=\frac{3\,(R')^2}{R^4}-\frac{3\,R''}{R^3},
&&R_{11}=\frac{2}{R^2}+\frac{(R')^2}{R^4}+\frac{R''}{R^3},
\quad\\
\vspace{-1ex}
&&&\mytag{6.16}\\
\vspace{-1ex}
&\hskip -2em
R_{22}=\frac{2}{R^2}+\frac{(R')^2}{R^4}+\frac{R''}{R^3},
&&R_{33}=\frac{2}{R^2}+\frac{(R')^2}{R^4}+\frac{R''}{R^3}.
\quad
\endxalignat
$$
And finally, using \mythetag{6.16}, we find the scalar curvature:
$$
R_{\,\text{scalar}}=-\frac{6}{R^2}-\frac{6\,R''}{R^3}.
\mytag{6.17}
$$
This formula \mythetag{6.17} coincides with the formula for 
the scalar curvature derived in \S\,112 of the book \mycite{1}.
\par
     2. {\bf Stereographic coordinates} in projection from the 
{\bf South Pole}. We denote these coordinates with $y^0,\,y^1,
\,y^2,\,y^3$. Their domain is the whole sphere $S^3$ 
except for the South Pole. The frames \mythetag{4.13} 
are defined and smooth in this domain. This frame pair is linked 
to the coordinates $y^0,\,y^1,\,y^2,\,y^3$. The commutation 
coefficients for the frame $\bold Y_0,\,\bold Y,\,\bold Y_2,
\,\bold Y_3$ in this pair are given by the formulas
\mythetag{1.19}. Other formulas \mythetag{6.10}, \mythetag{6.11}, 
\mythetag{6.12}, \mythetag{6.13}, \mythetag{6.14}, \mythetag{6.15}, 
\mythetag{6.16}, and \mythetag{6.17} are valid in this case upon
changing $x^0,\,x^1,\,x^2,\,x^3$ for $y^0,\,y^1,\,y^2,\,y^3$ in
them.\par
     3. {\bf Spherical coordinates}. We denote these coordinates 
with $\eta^0,\,\eta^1,\,\eta^2,\,\eta^3$ and also use special 
notations \mythetag{5.1} and \mythetag{5.2} for separate coordinates.
The domain of spherical coordinates is the whole sphere $S^3$ except 
for both poles North and South. The frame pair \mythetag{5.21} is used
for spherical coordinates. The frames of this pair are defined and
smooth in the domain of the spherical coordinates. The commutation 
coefficients for the frame $\bold E_0,\,\bold E,\,\bold E_2,
\,\bold E_3$ are given by the formulas \mythetag{5.9}. Substituting
them into \mythetag{6.5} we find the $\Gamma$-components of the
metric connection in the frame $\bold E_0,\,\bold E,\,\bold E_2,
\,\bold E_3$. Here is the list of nonzero ones of them:
$$
\xalignat 3
&\hskip -2em
\Gamma^0_{11}=\frac{R'}{R^2},
&&\Gamma^0_{22}=\frac{R'}{R^2},
&&\Gamma^0_{33}=\frac{R'}{R^2},
\quad\\
\vspace{1ex}
&\hskip -2em
\Gamma^1_{10}=\frac{R'}{R^2},
&&\Gamma^2_{20}=\frac{R'}{R^2},
&&\Gamma^3_{30}=\frac{R'}{R^2},
\quad\\
\vspace{-1ex}
&&&\mytag{6.18}\\
\vspace{-1ex}
&\hskip -2em
\Gamma^1_{22}=-\frac{\cos\chi}{R\,\sin\chi},
&&\Gamma^2_{33}=-\frac{\cos\theta}{R\,\sin\chi\sin\theta},
&&\Gamma^2_{21}=\frac{\cos\chi}{R\,\sin\chi},
\quad\\
\vspace{1ex}
&\hskip -2em
\Gamma^1_{33}=-\frac{\cos\chi}{R\,\sin\chi},
&&\Gamma^3_{32}=\frac{\cos\theta}{R\,\sin\chi\sin\theta},
&&\Gamma^3_{31}=\frac{\cos\chi}{R\,\sin\chi}.
\endxalignat
$$
Now we substitute \mythetag{6.18} into \mythetag{6.6} in order to
get the $\Alpha$-components  of the metric connection. Below is the 
list of nonzero ones of these components:
$$
\allowdisplaybreaks
\xalignat 2
&\hskip -2em
\Alpha^1_{12}=\frac{R'}{2\,R^2},
&&\Alpha^2_{11}=\frac{R'}{2\,R^2},
\quad\\
\vspace{1ex}
&\hskip -2em
\Alpha^3_{14}=-\frac{R'}{2\,R^2},
&&\Alpha^4_{13}=-\frac{R'}{2\,R^2},
\quad\\
\vspace{1ex}
&\hskip -2em
\Alpha^1_{21}=-\frac{i\,\cos\chi}{2\,R\,\sin\chi},
&&\Alpha^2_{22}=\frac{i\,\cos\chi}{2\,R\,\sin\chi},
\quad\\
\vspace{1ex}
&\hskip -2em
\Alpha^3_{23}=-\frac{i\,\cos\chi}{2\,R\,\sin\chi},
&&\Alpha^4_{24}=\frac{i\,\cos\chi}{2\,R\,\sin\chi},
\quad\\
\vspace{-1ex}
&&&\mytag{6.19}\\
\vspace{-1ex}
&\hskip -2em
\Alpha^1_{22}=-\frac{i\,R'}{2\,R^2},
&&\Alpha^2_{21}=\frac{i\,R'}{2\,R^2},
\quad\\
\vspace{1ex}
&\hskip -2em
\Alpha^3_{24}=\frac{i\,R'}{2\,R^2},
&&\Alpha^4_{23}=-\frac{i\,R'}{2\,R^2},
\quad\\
\vspace{1ex}
&\hskip -2em
\Alpha^1_{31}=\frac{R'}{2\,R^2},
&&\Alpha^2_{32}=-\frac{R'}{2\,R^2},
\quad\\
\vspace{1ex}
&\hskip -2em
\Alpha^3_{33}=-\frac{R'}{2\,R^2},
&&\Alpha^4_{34}=\frac{R'}{2\,R^2},
\quad\\
\endxalignat
$$
$$
\allowdisplaybreaks
\align
&\hskip -2em
\Alpha^1_{32}=\frac{\cos\chi}{2\,R\,\sin\chi}
-\frac{i\,\cos\theta}{2\,R\,\sin\chi\,\sin\theta},\\
\vspace{1ex}
&\hskip -2em
\Alpha^2_{31}=-\frac{\cos\chi}{2\,R\,\sin\chi}
-\frac{i\,\cos\theta}{2\,R\,\sin\chi\,\sin\theta},\\
\vspace{-1ex}
&\mytag{6.20}\\
\vspace{-1ex}
&\hskip -2em
\Alpha^3_{34}=\frac{\cos\chi}{2\,R\,\sin\chi}
-\frac{i\,\cos\theta}{2\,R\,\sin\chi\,\sin\theta},\\
\vspace{1ex}
&\hskip -2em
\Alpha^4_{33}=-\frac{\cos\chi}{2\,R\,\sin\chi}
-\frac{i\,\cos\theta}{2\,R\,\sin\chi\,\sin\theta}.
\endalign
$$
The formulas \mythetag{6.19} and \mythetag{6.20} are analogs of
\mythetag{6.11} and \mythetag{6.12}
Using the $\Gamma$-components of the metric connection \mythetag{6.18}
and applying the formula \mythetag{6.7} to them, we derive the explicit
formulas for the components of the Riemannian curvature tensor. Here is 
the list of nonzero components of $\bold R$:
$$
\xalignat 2
&\hskip -1ex
R^0_{101}=-R^0_{110}=-\frac{(R')^2}{R^4}+\frac{R''}{R^3},
&&R^1_{001}=-R^1_{010}=-\frac{(R')^2}{R^4}+\frac{R''}{R^3},
\qquad\quad\\
\vspace{1ex}
&\hskip -1ex
R^0_{202}=-R^0_{220}=-\frac{(R')^2}{R^4}+\frac{R''}{R^3},
&&R^2_{002}=-R^2_{020}=-\frac{(R')^2}{R^4}+\frac{R''}{R^3},
\qquad\quad\\
\vspace{1ex}
&\hskip -1ex
R^0_{303}=-R^0_{330}=-\frac{(R')^2}{R^4}+\frac{R''}{R^3},
&&R^3_{003}=-R^3_{030}=-\frac{(R')^2}{R^4}+\frac{R''}{R^3},
\qquad\quad\\
\vspace{-1ex}
&&&\qquad\quad
\mytag{6.21}\\
\vspace{-1ex}
&\hskip -1ex
R^1_{212}=-R^1_{221}=\frac{1}{R^2}+\frac{(R')^2}{R^4},
&&R^2_{112}=-R^2_{121}=-\frac{1}{R^2}-\frac{(R')^2}{R^4},
\qquad\quad\\
\vspace{1ex}
&\hskip -1ex
R^2_{323}=-R^2_{333}=\frac{1}{R^2}+\frac{(R')^2}{R^4},
&&R^3_{223}=-R^3_{232}=-\frac{1}{R^2}-\frac{(R')^2}{R^4},
\qquad\quad\\
\vspace{1ex}
&\hskip -1ex
R^3_{131}=-R^3_{113}=\frac{1}{R^2}+\frac{(R')^2}{R^4},
&&R^1_{331}=-R^1_{313}=-\frac{1}{R^2}-\frac{(R')^2}{R^4}.
\qquad\quad
\endxalignat
$$
Comparing \mythetag{6.21} with \mythetag{6.13}, we find that
these formulas are identical. Though written for two different 
frames $\bold X_0,\,\bold X_1,\,\bold X_2,\,\bold X_3$ and
$\bold E_0,\,\bold E_1,\,\bold E_2,\,\bold E_3$, the components
of the curvature tensor $\bold R$ do coincide. Then the same is 
true for the components of the Ricci tensor and the scalar 
curvature, i\.\,e\. they are given by the formulas \mythetag{6.16} 
and \mythetag{6.17} in spherical coordinates\footnotemark. Applying 
the formula \mythetag{6.9}, we conclude that the components of 
the spinor curvature tensor $\eufb R$ in spherical coordinates 
should coincide with those calculated for stereographic coordinates. 
They are given by the formulas \mythetag{6.14} and \mythetag{6.15}.
\footnotetext{\,We should warn here that typically the spherical 
coordinates are equipped with their canonical holonomic frame 
\mythetag{5.6}. We equip them with the non-holonomic frame 
\mythetag{5.7}.}

\enddocument
\end